\documentclass[11pt,letterpaper]{article}  
\usepackage{epsfig}
\usepackage{amsfonts}
\usepackage{amssymb}
\usepackage{amsmath}
\usepackage{float}
\def\qed{$\qquad \Box$}

\def\n1LXWo{{L^2_{X_0(0);W_0}(0,T;\mathbb{R}^{n_1})}}

\newtheorem{thm}{Theorem}[section]
\newtheorem{prop}[thm]{Proposition}
\newtheorem{lem}[thm]{Lemma}
\newtheorem{rem}[thm]{Remark}

\makeatletter \@addtoreset{equation}{section}

\renewcommand{\theequation}{\thesection.\arabic{equation}}
\newcommand{\beq}[1]{\begin{equation} \label{#1}}
\newcommand{\eeq}{\end{equation}}
\newcommand{\bed}{\begin{displaymath}}
\newcommand{\eed}{\end{displaymath}}

\newcommand{\bea}{\bed\begin{array}{rl}}
\newcommand{\eea}{\end{array}\eed}

\newcommand{\barray}{\begin{array}{ll}}
\newcommand{\earray}{\end{array}}

\newcommand{\bi}{\begin{itemize}}
\newcommand{\ei}{\end{itemize}}
\usepackage[normalem]{ulem}

\begin{document}

\vspace{-0.3cm}
\title{\LARGE\bf  Linear-Quadratic Mean Field Social Optimization\\ with a Major Player 
}

\vspace{1cm}
\author{Minyi Huang \thanks{M. Huang is with the School of
Mathematics and Statistics, Carleton University, Ottawa, ON K1S 5B6,
Canada (mhuang@math.carleton.ca). This work was supported in part  by Natural Sciences and Engineering
Research Council (NSERC) of Canada under a Discovery Grant and a
Discovery Accelerator Supplements Program.  Please address all correspondence to this author.}\quad  and Son Luu Nguyen
\thanks{S.L. Nguyen is with Department of Mathematics, University of Puerto Rico, San Juan, PR 00936 (sonluu.nguyen@upr.edu).}
     }


\maketitle

\begin{abstract}
This paper considers a linear-quadratic (LQ) mean field control problem involving a major player and a large number of minor players, where the dynamics and costs depend on random parameters.
 The objective is to optimize a social cost as a weighted sum of the individual costs  under decentralized information. We apply the  person-by-person optimality principle in team decision theory to the finite population model  to construct two limiting variational problems whose solutions, subject to the requirement of consistent mean field approximations, yield a system of forward-backward stochastic differential equations (FBSDEs). We show the existence and uniqueness of a solution to the FBSDEs and  obtain  decentralized strategies nearly achieving social optimality in the original large but finite population model.
\end{abstract}

\begin{center}
\begin{minipage}{15.3cm}
\end{minipage}
\end{center}

\begin{center}
\begin{minipage}{15.3cm}
{\bf Abbreviated title:} Mean field social optimization
\end{minipage}
\end{center}




\section{Introduction}

Mean field dynamic decision problems have been extensively studied in the literature \cite{BFY13,CHM17,CD18, HCM07,HCM12, HMC06, LL07, LZ08,  TB09},
and  a central goal  is to obtain decentralized strategies based on
limited information for individual agents. In a noncooperative game theoretic context,  decentralized solutions are developed in \cite{HCM07,HMC06} by applying consistent
mean field approximations.

In a basic mean field decision model, all players (or agents) have
comparably
small influence  and may be called peers. A modified modeling framework is to introduce one or a few  major players interacting with a large number of minor players. Traditionally, models differentiating the strength of players have
been  studied in cooperative game theory, and they are customarily called mixed
games  with the players according called mixed players \cite{Ha73}; such literature only dealt with static models. The  work \cite{H10} investigates an LQ mean field  game  involving a major player. The consideration of major and minor players in mean field control has attracted considerable interest addressing different nonlinear modeling aspects \cite{BLP14,CZ16,NC13}.
See \cite{BCY15,MB15} for extension to hierarchical games, and \cite{K17} for
the analysis of evolutionary
games and their deterministic mean field limit under a principle agent.

On the other hand, cooperation in dynamic multi-agent decision problems
is traditionally a well studied subject. For general cooperative  differential games using various optimality notions, see \cite{SM76,THL86,YP06}.
Naturally, cooperative decision making in mean field models is  of interest, especially from the point of view of addressing complexity
\cite{HCM12,TB09}. Such decision problems may be referred to as mean field teams for which the decision makers will also be called players or agents.  The work \cite{HCM12} introduced an LQ social optimization problem  where all the agents cooperatively minimize a social cost as the sum of their individual costs, and it  shows that the consistency based approach in mean field games may be  extended to this model by combining with a person-by-person optimality principle in team decision theory \cite{H80,WS00}.  The central result is the so-called social optimality theorem which states that the optimality loss of the obtained decentralized strategies becomes negligible when the population size goes to infinity \cite{HCM12}.
The social optimum may be regarded as a specific Pareto optimum for the constituent agents.
 A mean field team is studied in \cite{WZ17} where
a Markov jump parameter appears as a common source of
randomness for all agents.
An LQ mean field team is formulated in  \cite{AM15}
by assuming mean field (i.e. the average state of the
population) sharing for a given population size $N$, which gives an optimal control problem with special partial state information. In a mixed player setting, \cite{BLP14} considers a nonlinear diffusion model and
assumes that all minor players act as a team to minimize a common cost
against the major player. Optimal control of McKean--Vlasov dynamics is analyzed in \cite{L17} and under some conditions it is shown that the optimal solution may be interpreted as the limit of the social optimum solution of
$N$-players as $N\to \infty$.  Cooperative mean field control has applications in economic theory \cite{NM18}, collective motion control \cite{ACFK17,PRT15}, and power grids \cite{CBM17}.
Furthermore,  social optima are useful for studying efficiency of mean field games by providing a performance benchmark \cite{BT13,HCM03}.

For mean field teams with mixed players,
the analysis in an LQ framework has been formulated in our earlier work \cite{HN11}, where  partial analysis was presented by applying a state space  augmentation technique  to characterize the dynamics of the random mean
field evolution. Later,
 \cite{HN16} re-examined the problem  by applying the person-by-person optimality principle adopted for the  peer model in \cite{HCM12}. This paper further generalizes the model by including coupling in dynamics and random coefficients while  \cite{HN11} only considers cost coupling and deterministic coefficients. Specifically, the model parameters now depend on the Brownian motion of the major player. This suggests that  the major player serves as a  common source of randomness for all players, which has connections with  mean field games with common noise \cite{BFY13,CDLL15,CD18}. In fact,
 the stochastic control literature  \cite{B78,P92} has considered a similar randomness structure where the system  coefficients depend on a smaller filtration, and such modeling has applications in finance \cite{KZ00,L04}.

 As in \cite{HN16}, our solution is to extend the person-by-person optimality argument
of \cite{HCM12} to the current setting to
deal with random mean field approximations due to the presence of the major player,
and we  solve two variational problems resulting from  the major-minor
player interactions. The linear backward stochastic differential equation (BSDE) \cite{B76,MY99} technique adopted in this paper can treat the random mean field and coefficients in a unified manner.
 As it turns out,  the consideration of the coupling in dynamics will necessitate
 delicate handling of a two-scale variational problem for the minor player.
Note that for the person-by-person optimality principle
only one player has control perturbation. This  feature is similar to mean field games where the equilibrium is tested by unilateral strategy changes. However, our performance characterization of social optimality must allow simultaneous control variations. The optimal control nature of our problem shares some similarity with mean field type optimal control \cite{ELN13,Y13}. However, the later involves only a single decision maker which directly controls the state mean.

Throughout this paper, we use  $(\Omega, {\cal F}, \{{\cal F}_t\}_{t\ge 0}, P)$ to denote an underlying filtered probability space.
Let $S^n$ be the Euclidean space of $n\times n$ real and symmetric matrices,  $S_+^n$ its subset of
positive semi-definite matrices, and $I_k$ the $k\times k$ identity matrix.
 The Banach space $L_{\cal F}^2(0, T; \mathbb{R}^k)$ consists of all
 $\mathbb{R}^k$-valued ${\cal F}_t$-adapted
  square integrable processes $\{v(t), 0\le t\le T\}$
   with the norm $\|v\|_{L^2_{\cal F}}=( E\int_0^T|v(t)|^2 dt)^{1/2}$.
The Banach space $L_{\cal F}^\infty(0, T; \mathbb{R}^k)$ consists of
    all $\mathbb{R}^k$-valued ${\cal F}_t$-adapted
  essentially bounded processes $\{v(t), 0\le t\le T\}$
   with the norm $\|v\|_{L^\infty_{\cal F}}=\mbox{ess} \sup_{t, \omega} |v(t)|$.
   We can similarly define such spaces with other choices of the filtration
   and the Euclidean space.
Given a symmetric matrix $M\geq 0$, the quadratic form $z^TMz$ may
be denoted as $|z|_M^2$. For a matrix $Z$, $Z_i^{\rm col}$ stands for the $i$th column of $Z$.
Some variables (such as $X_0^\star(t)$, $u_i^\star(t)$) with a
superscript of star are used for limiting models afer taking mean field approximations.
Let $\{{\cal F}_t^W, t\ge 0\}$ be the filtration by a Brownian motion  $\{W(t), t\ge 0\}$. We use $C$ (or $C_1$, etc.) to denote a generic constant which is independent of the population size $N+1$ and may change from place to place.


The organization of the paper is as follows. Section \ref{sec:mod} formulates the social optimization problem with a major player. Sections \ref{sec:maj} and \ref{sec:min} introduce two variational problems with random parameters  for the major player and a representative minor player, respectively. The existence and uniqueness of the  mean field social optimum solution is presented in section \ref{sec:mfsol}. An asymptotic social optimality theorem is established in section \ref{sec:asy}. Section \ref{sec:con} concludes the paper.

 \section{The Mean Field Social Optimization Model}

\label{sec:mod}

Consider the LQ mean field decision model  with a major player ${\cal A}_0$  and   minor players $\{{\cal A}_i, 1\leq i\leq N\}$. 
 At time $t\geq 0$, the states of ${\cal A}_0$  and ${\cal A}_i$ are, respectively, denoted by $X_0^N(t)$ and $X_i^N(t)$, $1\leq i\leq N$.
The dynamics of the $N+1$ players are given by a system of linear stochastic differential equations (SDEs):
  \begin{align}
dX_0^N(t) =& \Big[A_0(t) X_0^N(t)+B_0(t) u_0^N(t) +F_0(t)X^{(N)}(t) \Big]dt   +D_0(t) dW_0(t), \qquad   \label{maj}\\
dX_i^N(t)=&  \Big[A(t) X_i^N(t) +B(t)u_i^N(t)+F(t)X^{(N)}(t)+G(t)X_0^N(t) \Big] dt
       +  D(t)dW_i(t), \label{min}
        \quad 1\leq i\leq N,
\end{align}
where $X^{(N)}(t)=(1/N)\sum_{i=1}^N X_i^N(t)$ is the coupling term.
The states $X_0^N$, $X_i^N$ and controls $u_0^N$, $u_i^N$ are, respectively, $n$ and $n_1$ dimensional vectors. The initial states $X_0^N(0)= z_{0}$, $X_i^N(0)=x_{i0}^N$, $1\leq i\leq N$,  are deterministic. The coefficients in the dynamics are random.      The noise processes $W_0$, $W_i$ are  $n_2$ dimensional independent standard Brownian motions adapted to ${\cal F}_t$.
  We  choose ${\cal F}_t$  as the $\sigma$-algebra
  $ {\cal F}_t^{W_{\centerdot}}
  \mathrel{\mathop:}= \sigma(W_j(\tau), 0\leq j\leq N,\tau\leq t)$.
Denote  $W_0=[W_{01},\ldots, W_{0n_2}]^T$,
   ${\cal F}_t^{W_0}\mathrel{\mathop:}=\sigma (W_0(\tau), \tau \le t)$,
   and ${\cal F}_t^{W_0,W_i}\mathrel{\mathop:}=\sigma (W_0(\tau),W_i(\tau), \tau \le t)$.


For $0\leq j\leq  N$, denote $u_{-j}^N=\big(u_0^N,\ldots, u_{j-1}^N,u_{j+1}^N,\ldots, u_N^N\big)$.
The cost  for ${\cal A}_0$ is given by
\begin{align}
J_0(u_0^N, u_{-0}^N)
=\ &  E\int_0^T\Big\{ \big|X_0^N(t)- \Psi_0(X^{(N)}(t))\big|_{Q_0(t)}^2
 + (u_0^N(t))^T R_0(t) u_0^N(t) \Big\}dt \nonumber \\
 &\quad +E |X_0^N(T)- H_{0,f} X^{(N)}(T)-\eta_{0f}|_{Q_{0f}}^2, \label{lc}
\end{align}
where
   $\Psi_0(X^{(N)}(t))=H_0(t)X^{(N)}(t)+\eta_0(t)$.
The cost for ${\cal A}_i$, $1\leq i\leq N$, is given by
\begin{align}
 J_i(u_i^N, u_{-i}^N)  =\ &
 E\int_0^T \Big\{\big|X_i^N(t)-\Psi(X_0^N(t),X^{(N)}(t))\big|_{Q(t)}^2
  + (u_i^N(t))^T R(t) u_i^N (t)\Big\}dt \nonumber \\
 & \quad+ E|X_i^N(T)-H_{1f} X_0^N(T)-  H_{2f}X^{(N)}(T)-\eta_{f} |_{Q_f}^2,  \label{fc}
\end{align}
where $\Psi(X_0^N(t),X^{(N)}(t))= H_1(t) X_0^N(t)+  H_2(t)X^{(N)}(t)+\eta(t)$. The terms $H_1(t)X_0^N(t)$ and $H_{1f} X_0^N(T)$ indicate the strong influence of the major player. Also, the parameters in the two costs are random.

Below we list the stochastic parameter  processes
     \begin{align}
&\{A_0(t),\ B_0(t),\ F_0(t),\  D_0(t),\ A(t),\ B(t),\ F(t),\ G(t), \ D(t), 0\le t\le T\},  \label{AD} \\
&\{ H_0(t),\ H_1(t),\  H_2(t),\  Q_0(t),\  Q(t) ,\ R_0(t),\  R(t),\ \eta_0(t),\
    \eta(t), 0\le t\le T\}.   \label{Het}
    \end{align}

 We introduce the standing assumptions for this paper.

(A1) We have
\begin{align*}
&A_0, F_0, A, F, G , H_0, H_1, H_2 \in L^\infty_{{\cal F}^{W_0}} (0,T; \mathbb{R}^{n\times n}), \qquad\\
&B_0, B\in   L^\infty_{{\cal F}^{W_0}} (0,T; \mathbb{R}^{n\times n_1}), \quad
D_0, D\in  L^2_{{\cal F}^{W_0}} (0,T; \mathbb{R}^{n\times n_2}), \\
& \eta_0, \eta \in L^2_{{\cal F}^{W_0}} (0,T; \mathbb{R}^n),
\end{align*}
    and
\begin{align*}
& Q_0, Q\in L^\infty_{{\cal F}^{W_0}} (0,T; S^n), \quad Q_0(t)\in S_+^n,\ Q(t) \in S_+^n, \quad \forall t\in [0,T], \\
& R_0, R\in L^\infty_{{\cal F}^{W_0}} (0,T; S^{n_1}), \quad R_0(t)\ge c_1 I_{n_1}, \quad R(t)\ge c_1 I_{n_1}, \qquad \forall t\in[0,T],
\end{align*}
where $c_1>0$ is a fixed deterministic constant.

(A2)
The terminal cost parameters
\begin{align}
 H_{0f}, \ Q_{0f}, \ H_{1f},\ H_{2f},\ Q_f, \label{Hetf}
\end{align}
are ${\cal F}_T^{W_0}$-measurable and essentially bounded, and $ Q_{0f},  Q_f$ are $S_+^n$-valued. $ \eta_{0f}$ and $  \eta_f $ are
 ${\cal F}_T^{ W_0}$-measurable and square integrable.




(A3)
There exists a constant $c_{2}>0$ independent of $N$ such that
$\sup_{j\ge 0}| x_{j0}^N|\le c_2$ for the initial states, and $\lim_{N\to \infty} x_0^{(N)}=m_0$,
where $x_0^{(N)}=\tfrac{1}{N} \sum_{i=1}^N x_{i0}^N$.

By (A1)--(A2),
there exists a fixed constant $c_3$ such that
$$
\mbox{ess}\sup_{t, \omega}|\psi(t)|\le c_3, \quad \mbox{ess} \sup_\omega | \psi_f|\le c_3,
$$
where $\psi(t)$ (resp., $\psi_f$) stands for
any entry in \eqref{AD}--\eqref{Het} (resp., \eqref{Hetf}).

For the rest of the paper, for a stochastic process $\{Z(t), 0\le t\le T\}$ appearing in various equations and equalities, we may write $Z$ for $Z(t)$ by suppressing the time variable $t$ for which the interpretation should be clear from the context. For instance, we often drop $t$ in $A_0(t)$, $B_0(t)$, $X_i^N(t)$, $Q(t)$, etc.

Throughout the paper, we denote $Y^{(N)}= \frac{1}{N}\sum_{i=1}^N
Y_i$, and $Y_{-i}^{(N)}= \frac{1}{N}\sum_{j\ne i}^N
Y_j$ for $N$ vectors $Y_1, \ldots, Y_N$.

\subsection{The mean field social optimization problem}

For the mean field  social optimization problem, we attempt  to minimize the following social cost
\begin{align}
J_{\rm soc}^{(N)}(u)=J_0+\frac{\lambda}{N}\sum_{k=1}^N J_k,
\end{align}
where $ u^N=(u_0^N,u_1^N, \ldots, u_N^N) $ and $\lambda>0$. It is necessary to
introduce the scaling factor $\lambda/N$ in order to obtain a well defined
limiting problem when $N$ tends to infinity.
In view of the dynamics and costs of the $N+1$ players, $J_0$ and $J_i$, $i\geq 1$,
 are generally of the same order of magnitude. If $\lambda/N$ were replaced by 1, the limiting
control problem would be too insensitive to the performance of the major player and
become inappropriate.

For the model of $N+1$ players, let the optimal control be denoted by
\begin{align}
\check u^N= \big(\check u_0^N, \check u_1^N, \ldots, \check u_N^N\big), \label{teamop}
\end{align}
where each $\check u_j$  belongs to $L_{\cal F}^2(0, T; \mathbb{R}^{n_1})$.
Since the optimal control problem minimizing $J_{\rm soc}^{(N)}$ is a strictly convex optimization
problem with $J_{\rm soc}^{(N)} \to \infty $ as $\|u^N\|_{L^2_{\cal F}}\to \infty$, such $\check u^N$ exists and is unique.
However, this  solution is not what  we desire to obtain since each player needs centralized information. Instead, it will serve as a starting point for designing decentralized strategies.


\section{The Major Player's Variational Problem}

\label{sec:maj}

Consider the variation $\tilde u_0^{N} \in L^2_{\cal F}(0, T; \mathbb{R}^{n_1}) $ and $u_0^N=\check u_0^N+\tilde u_0^N$. Let the state processes $(\check X_j^{N})_{j=0}^N$ correspond to
$(\check u_j^{N})_{j=0}^N$, and
$(X_j^N)_{j=0}^N$ correspond to
$(\check u_0^{N} +\tilde u_0^{N}, \check u_1^{N}, \ldots, \check u_N^{N} )$. Write $X_j^N=\check X_j^{N} +\tilde X_j^{N}$ for $0\le j\le N$, where $\tilde X_j^{N}$ is the state variation of player $\mathcal{A}_j$.
Then
\begin{align*}
&d X_0^{N}(t) = (A_0 X_0^{N} + F_0 X^{(N)} +B_0 u_0^{N} )dt+D_0dW_0(t), \\
& d X^{(N)}(t) = [(A+F) X^{(N)} +B\check u^{(N)}+G X_0^{N}] dt+\frac{D}{N} \sum_{i=1}^N dW_i(t),
\end{align*}
and
\begin{align*}
&d \tilde X_0^{N}(t) = (A_0 \tilde X_0^{N} + F_0 \tilde X^{(N)} +B_0 \tilde u_0^{N} )dt, \\
& d\tilde X^{(N)}(t) = [(A+F) \tilde X^{(N)} +G \tilde X_0^{N}] dt,
\end{align*}
where $\tilde X_0^N (0)= \tilde X^{(N)}(0)=0$.
Note that we have followed the convention of dropping  the time  variable $t$ in various places.
It can be checked that $\tilde X_i^N= \tilde X^{(N)}$ on $[0,T]$ for all $1\le i\le N$.
Denote
\begin{align*}
\delta L_0^N(t)=&\Big\{ [\check X_0^{N} - \big(H_0 \check X^{(N)}+\eta_0\big)]^T Q_0 \big(\tilde X_0^{N} -H_0 \tilde X^{(N)}\big)+
 \big(\check u_0^{N}\big)^T R_0 \tilde u_0^{N} \\
&+\lambda [(I- H_2)\check  X^{(N)}- H_1 \check X_0^{N} -\eta]^T Q[(I- H_2)
\tilde X^{(N)} -H_1 \tilde X_0^{N}]\Big\}(t),
\end{align*}
and
\begin{align*}
\delta L_{0f}^N=& \Big\{[\check X_0^{N} - \big(H_{0f} \check X^{(N)}+\eta_{0f}\big)]^T Q_{0f} \big(\tilde X_0^{N} -H_{0f} \tilde X^{(N)}\big) \\
&+\lambda [(I- H_{2f})\check  X^{(N)}- H_{1f} \check X_0^{N} -\eta_f]^T Q_f[(I- H_{2f}) \tilde X^{(N)} -H_{1f} \tilde X_0^{N}]\Big\}(T).
\end{align*}

The first variation of the social cost is given by
\begin{align*}
\delta J_0 + \frac{\lambda}{N} \sum_{i=1}^N \delta J_i =
2E\int_0^T \delta L_0^N (t) dt  +2E\delta L_{0f}^N  ,
\end{align*}
which is a linear functional of $\tilde u_0^{N}$. We have the first order variational condition:

\begin{lem}\label{lemma:pbp0}
We have
\begin{align*}
E\int_0^T \delta L_0^N(t) dt+E\delta L_{0f}^N  =0,\quad \forall\  \tilde u_0^{N} \in L^2_{\cal F}(0, T; \mathbb{R}^{n_1}).
\end{align*}
\end{lem}

{\it Proof.} We prove by using the so called  person-by-person optimality
principle \cite{H80}. Take a constant $\epsilon$ and let $\tilde u_0^N$ be fixed. Then consider the control
$(\check u_0^{N} +\epsilon \tilde u_0^{N}, \check u_1^{N}, \ldots, \check u_N^{N} )$ for the players. It follows that
$$
J_{\rm soc}^{(N)}(\check u_0^{N} +\epsilon \tilde u_0^{N}, \check u_1^{N}, \ldots, \check u_N^{N}) \ge J_{\rm soc}^{(N)}(\check u_0^{N} , \check u_1^{N}, \ldots, \check u_N^{N})
$$
for all $\epsilon$, and the lemma follows from elementary estimates of the left hand side after an expansion around $\check u_0$. \qed

\subsection{The limiting variational problem for the major player}
Consider the limiting model
\begin{align*}
dX_0^{\star}(t)&= (A_0 X_0^{\star}+B_0 u_0^{\star} +F_0 m )dt+D_0 dW_0(t), \\
dm(t) &= [(A+F) m +B\bar u  +G X_0^{\star}]dt,
\end{align*}
where  $X_0^\star(0)= X_0^N(0)$, $m(0)=m_0$,
 and $\bar u\in L_{{\cal F}^{W_0}}^2(0,T; \mathbb{R}^{n_1})$.  Here  $\bar u $ and $m$ are used to approximate
$\check u^{(N)}$ and $X^{(N)}$ for large $N$, respectively.
Note that each $\check u_j \in L^2_{\cal F}(0,T;\mathbb{R}^{n_1})$ is a  centralized control in that it depends on all  Brownian motions. However, as $N\to \infty$, we expect the randomness originated
 in  $ (W_1,
\ldots, W_N)$ will be averaged out. This has motivated the consideration of $\bar u \in L_{{\cal F}^{W_0}}^2(0,T; \mathbb{R}^{n_1})$.

For a particular control $\hat u_0^{\star}\in  L_{{\cal F}^{W_0}}^2(0,T; \mathbb{R}^{n_1})$,  let the associated state process be
\begin{align*}
d\hat X_0^{\star}(t)&= (A_0 \hat X_0^{\star}+B_0 \hat u_0^{\star} +F_0\hat m)dt+D_0 dW_0(t), \\
d\hat m(t)&= [(A+F) \hat m +B\bar u +G \hat X_0^{\star}] dt,
\end{align*}
where $\hat X_0^\star(0) =X_0^N(0)$ and $\hat m(0)=m_0.$
The state variations read
\begin{align*}
d \tilde X_0^{\star}(t)&=(A_0 \tilde X_0^{\star}  + F_0 \tilde m +B_0\tilde u_0^{\star})dt,  \\
d\tilde m(t)&=[(A+F)\tilde m+G\tilde X_0^{\star}] dt,
\end{align*}
where $\tilde X_0^\star(0)=0$, $\tilde m(0)=0$, and  $\tilde u_0^{\star} \in L_{{\cal F}^{W_0}}^2(0,T; \mathbb{R}^{n_1}) $ is the control variation.

 Denote
\begin{align*}
K_0(t)&=-Q_0H_0-\lambda H_1^TQ(I- H_2),\\
M_0(t)&=Q_0+\lambda H_1^TQH_1,\\
M(t)&=H_0^TQ_0H_0+\lambda(I- H_2)^TQ(I- H_2),\\
\nu_0(t)&=\lambda H_1^TQ\eta-Q_0\eta_0,\\
\nu(t)&=H_0^TQ_0\eta_0+\lambda  H_2^TQ\eta-\lambda Q\eta,\\
R_\lambda(t)&= \lambda R,
\end{align*}
where the time variable $t$ in various places of the right hand sides is suppressed.
For the terminal cost, similarly define
\begin{align*}
K_{0f}&=-Q_{0f}H_{0f}-\lambda H_{1f}^TQ_f(I- H_{2f}),\\
M_{0f}&=Q_{0f}+\lambda H_{1f}^TQ_fH_{1f},\\
M_f&=H_{0f}^TQ_{0f}H_{0f}+\lambda(I- H_{2f})^TQ_f(I- H_{2f}),\\
\nu_{0f}&=\lambda H_{1f}^TQ_f\eta_f-Q_{0f}\eta_{0f},\\
\nu_f&=H_{0f}^TQ_{0f}\eta_{0f}+\lambda  H_{2f}^TQ_f\eta_f-\lambda Q_f\eta_f.
\end{align*}

Define
\begin{align*}
\delta L_0^{\star}(t)=\ &\Big\{ \big[\hat X_0^{\star}-
\big(H_0 \hat m+\eta_0\big)\big]^T Q_0(\tilde X_0^{\star}-H_0 \tilde m)+
 (\hat u_0^{\star})^T R_0 \tilde u_0^{\star} \\
&+\lambda [(I- H_2)\hat  m- H_1\hat X_0^{\star} -\eta]^T Q[\big(I- H_2\big)
\tilde m -H_1 \tilde X_0^{\star}] \Big\}(t)  \\
=\ & \Big\{ (\tilde X_0^\star)^T ( M_0 \hat X_0^\star + K_0 \hat m + \nu_0) +
\tilde m^T ( K_0^T \hat X_0^\star + M \hat m +\nu       )+(\hat u_0^{\star})^T R_0 \tilde u_0^{\star} \Big\}(t),
\end{align*}
$$
\delta L_{0f}^\star= \{(\tilde X_0^\star)^T ( M_{0f} \hat X_0^\star + K_{0f} \hat m + \nu_{0f}) +
\tilde m^T ( K_{0f}^T \hat X_0^\star + M_f \hat m +\nu_f       )\}(T),
$$
which are intended to approximate $\delta L_0^N(t)$ and  $\delta L_{0f}^N$,
respectively.

{\bf Variational Problem (I)} VP--(I):
Find $ \hat u_0^{\star}\in L_{{\cal F}^{W_0}}^2(0,T; \mathbb{R}^{n_1})    $  such that
\begin{align} \label{vp0}
E\int_0^T \delta L_0^\star(t) dt +E\delta L_{0f}^\star  =0, \quad \forall \ \tilde u_0^{\star} \in L_{{\cal F}^{W_0}}^2(0,T; \mathbb{R}^{n_1})  .
\end{align}
We call $\hat u_0^\star$ or $(\hat u_0^{\star}, \hat X_0^{\star}, \hat m)$ a solution of { VP}--(I).

Suppose  $(\hat u_0^{\star}, \hat X_0^{\star}, \hat m)$ is a solution to
VP--(I).
We introduce the  backward stochastic differential equations (BSDEs)
\begin{align}
 dp_0(t)
 &=(M_0\hat X_0^\star+K_0 \hat m-A_0^Tp_0-G^Tp+\nu_0)dt+\xi_0dW_0(t),\label{eqnp0}\\
dp (t)
&=[K_0^T\hat X_0^\star+M\hat m-F_0^Tp_0-(A+F)^Tp+\nu]dt +\xi dW_0(t),\label{eqnp}
\end{align}
where
\begin{align}\label{ppT}
p_0(T)=-( M_{0f} \hat X_0^\star(T) + K_{0f} \hat m(T) + \nu_{0f}),\quad  p(T)= -( K_{0f}^T \hat X_0^\star(T) + M_f \hat m(T) +\nu_f       ) .
 \end{align}

\begin{lem}\label{lemma:urep}
If  $(\hat u_0^{\star}, \hat X_0^{\star}, \hat m)$  is a solution to {\em VP--(I)}, then \eqref{eqnp0}--\eqref{eqnp} has a unique solution $(p_0,p, \xi_0, \xi)$   in $L_{{\cal F}^{W_0}}^2(0,T; \mathbb{R}^{2n})  \times
L_{{\cal F}^{W_0}}^2(0,T; \mathbb{R}^{2n \times n_2})  $,   and
$\hat u_0^\star(t) = R_0^{-1}(t) B_0^T(t) p_0(t).$
\end{lem}

{\it Proof}. From the linear BSDEs, we can solve a unique solution
$(p_0, p, \xi_0, \xi)$.
 It follows from It$\hat{\mbox{o}}$'s formula that
\begin{align*}
&d[p_0^T(t)\tilde X_0^{\star}(t)+p^T(t)\tilde m(t)]\\
=\,\,&p_0^T(A_0 \tilde X_0^{\star}  + F_0 \tilde m +B_0\tilde
u_0^{\star})dt+p^T[(A+F)\tilde m+G\tilde X_0^{\star}] dt\\
&\quad+ (\tilde X_0^\star)^T (M_0\hat X_0^\star+K_0\hat m
-A_0^Tp_0-G^Tp+\nu_0)dt+(\tilde X_0^\star)^T\xi_0dW_0(t) \\
& \quad+\tilde m^T [K_0^T\hat X_0^\star+M\hat m-F_0^Tp_0-(A+F)^Tp+\nu]dt +\tilde m^T \xi dW_0(t).
\end{align*}
 Since $\tilde X_i^\star (0)=\tilde m(0)=0$,  this  implies
\begin{align}\label{pX0T}
   &  E\Big[p_0^T(T)\tilde X_0^{\star}(T)+p^T(T)\tilde m(T)\Big] \nonumber \\
   & = E\int_0^T[
 p_0^TB_0 \tilde u_0^{\star} + (\tilde X_0^\star)^T (M_0 \hat X_0^\star +K_0 \hat m+\nu_0)  + \tilde m^T (K_0^T \hat X_0^\star +M \hat m+\nu)   ] dt.
 \end{align}
It follows from \eqref{vp0} and \eqref{pX0T} that for any
$\tilde u_0^{\star}\in L_{{\cal F}^{W_0}}^2(0,T; \mathbb{R}^{n_1})$,
\begin{align} \label{ubp}
E\int_0^T (\tilde u_0^{\star})^T   \left(B_0^Tp_0 -   R_0 \hat u_0^\star  \right)
  dt=0.
\end{align}
The lemma follows. \qed

Given $\bar u\in L_{{\cal F}^{W_0}}^2(0,T; \mathbb{R}^{n_1})$,  denote the forward-backward stochastic differential equation (FBSDE)
\begin{equation}\label{clmaj}
\begin{cases}
d\hat X_0^{\star}(t)= (A_0 \hat X_0^{\star}+B_0R_0^{-1}B_0^Tp_0 +F_0\hat m)dt+D_0 dW_0(t), \\
d\hat m(t)= [(A+F) \hat m +B\bar u  +G \hat X_0^{\star}] dt,\\
dp_0(t)  =\big(M_0\hat X_0^\star+K_0\hat m-A_0^Tp_0-G^Tp+\nu_0\big)dt+\xi_0dW_0(t),\\
dp (t)  =[K_0^T\hat X_0^\star+M\hat m-F_0^Tp_0-(A+F)^Tp
+\nu]dt +\xi dW_0(t),
\end{cases}
\end{equation}
where $\hat X_0^\star(0)=X_0^N(0)$, $\hat m(0)=m_0$, $p_0(T)=-( M_{0f} \hat X_0^\star(T) + K_{0f} \hat m(T) + \nu_{0f})$,  $p(T)= -( K_{0f}^T \hat X_0^\star(T) + M_f \hat m(T)+\nu_f)$.

To analyze  \eqref{clmaj}, we introduce the notation:
\begin{align}
&{\bf X}_0=
\left[
\begin{matrix}
\hat X_0^\star\\
\hat m
\end{matrix}
\right],
\quad
{\bf Y}_0=
\left[
\begin{matrix}
p_0\\
p
\end{matrix}
\right],
\quad
Z_0=
\left[
\begin{matrix}
\xi_0\\
\xi
\end{matrix}
\right],\quad
{\mathbb A}_0=
\left[
\begin{matrix}
A_0&F_0\\
G&A+F
\end{matrix}
\right],
\quad
{\mathbb B}_0=
\left[
\begin{matrix}
B_0\\
0
\end{matrix}
\right],
\qquad
   \label{not1}    \\
& {\mathbb B}=
\left[
\begin{matrix}
0\\
B
\end{matrix}
\right],
\quad
{\mathbb D}_0=
\left[
\begin{matrix}
D_0\\
0
\end{matrix}
\right],\quad
{\mathbb Q}_0=
\left[
\begin{matrix}
M_0&K_0\\
K_0^T&M
\end{matrix}
\right],\quad
{\bf v}_0 =
\left[
\begin{matrix}
\nu_0\\
\nu
\end{matrix}
\right],
\quad \label{not2}  \\
&{\mathbb Q}_{0f}=
\left[
\begin{matrix}
M_{0f}&K_{0f}\\
K_{0f}^T&M_f
\end{matrix}
\right],\quad {\bf v}_{0f} =
\left[
\begin{matrix}
\nu_{0f}\\
\nu_f
\end{matrix}
\right]. \label{not3}
\end{align}
\begin{lem} \label{lemma:Q0}
 ${\mathbb Q}_0(t)$ and  ${\mathbb Q}_{0f}$  are  positive semi-definite for all $0\le t\le T$.
\end{lem}
{\it Proof.} Since $Q_0(t),Q(t)$ are symmetric and $Q_0(t),Q(t)\ge0$, we can write $Q_0(t)=U_0^T(t)U_0(t)$ and $Q(t)=U^T(t)U(t)$ for some
$\mathbb{R}^{n\times n}$-valued random matrices $U_0(t),U(t)$. Denote
$$
{\mathbb U}_0(t)=
\left[
\begin{matrix}
U_0&-U_0(t)H_0\\
0&0
\end{matrix}
\right](t),
\qquad
{\mathbb U}(t)=\sqrt\lambda
\left[
\begin{matrix}
UH_1&-U\big(I- H_2\big)\\
0&0
\end{matrix}
\right](t).
$$
It is clear that ${\mathbb Q}_0(t)$ is symmetric and
\begin{align*}
{\mathbb Q}_0(t)
&=\left[
\begin{matrix}
Q_0&-Q_0H_0\\
-H_0^TQ_0&H_0^TQ_0H_0
\end{matrix}
\right](t)
+\lambda
\left[
\begin{matrix}
H_1^TQH_1&-H_1^TQ\big(I- H_2\big)\\
-\big(I- H_2\big)^TQH_1&\big(I- H_2\big)^TQ\big(I- H_2\big)
\end{matrix}
\right](t) \\
&={\mathbb U}_0^T(t) {\mathbb U}_0(t)+{\mathbb U}^T(t) {\mathbb U}(t)\ge 0.
\end{align*}
The case of ${\mathbb Q}_{0f}$ can be similarly checked.
\qed

\begin{thm} \label{theorem:4equ}
i) For any $\bar u\in L_{{\cal F}^{W_0}}^2(0,T; \mathbb{R}^{n_1})$,
\eqref{clmaj} has a unique solution $(\hat X_0^\star,\hat m, p_0, p, \xi_0, \xi)$ in $ L_{{\cal F}^{W_0}}^2(0,T; \mathbb{R}^{4n})   \times L_{{\cal F}^{W_0}}^2 (0,T;\mathbb{R}^{2n\times n_2})$.

ii) {\em VP-(I)} has a unique solution given by
\begin{align}\label{optu}
\hat u_0^\star(t)=R_0^{-1}(t)B_0^T(t)p_0(t).
\end{align}
\end{thm}

{\it Proof.} i)
 We  rewrite \eqref{clmaj} in the form
$$
\begin{cases}
d{\bf X}_0(t)=\big({\mathbb A}_0{\bf X}_0+{\mathbb B}_0{ R}_0^{-1}{\mathbb B}_0^T{\bf Y}_0+ \mathbb{B}\bar { u}\big)(t)dt+{\mathbb D}_0(t)dW_0(t),\\
d{\bf Y}_0(t)=\big({\mathbb Q}_0{\bf X}_0-{\mathbb A}_0^T{\bf Y}_0+{\bf v}_0\big)(t)dt+{ Z_0}(t)dW_0(t),
\end{cases} $$
where ${\bf Y}_0(T)=-{\mathbb Q}_{0f} {\bf X}_0(T) - {\bf v}_{0f}$.
Under (A1)--(A2) and in view of Lemma  \ref{lemma:Q0}, we apply Lemma \ref{lemma:fbXYZ} to
obtain the existence and uniqueness of a solution.

ii) We solve
\eqref{clmaj} and choose $\hat u_0^\star$ by \eqref{optu}.
Such a control $\hat u_0^\star$ ensures \eqref{ubp} while  \eqref{pX0T} still holds; this further implies \eqref{vp0}. Hence, $\hat u_0^\star$ is a solution to VP--(I).

 On the other hand, if $(\hat u_0^\star, \hat X_0, \hat m)$ is a solution to VP--(I) such that \eqref{vp0} holds,
by Lemma \ref{lemma:urep}, \eqref{clmaj} holds and is uniquely solved, and uniqueness of $\hat u_0^\star $ follows from its representation in Lemma \ref{lemma:urep}.~\qed

\section{The Minor Player's Variational Problem}

\label{sec:min}


  Recall that the control $\check u$ yields state processes $
\check X_j^{N}$, $j=0, \ldots, N$.
Now consider the control
$(u_i^{N}, \check u_{-i}^N)$ for a fixed  $i\ge 1$.
Note that the state of player ${\cal A}_j$, $1\le j\ne i\le N$, is
affected even if only $\check u_i^N$ changes to $u_i^N
$. Under $(u_i^{N}, \check u_{-i}^N)$,
the state process of  player $\mathcal{A}_j$, $1\le j \ne i\le N$, is
$$
d  X_j^{N}(t) = ( A  X_j^{N} +B \check u_j^{N} +F  X^{(N)}+G  X_0^{N})dt  +D d W_j(t).
$$
Thus,
\begin{align*}
d  X_0^{N}(t)&= \Big( A_0  X_0^{N} + B_0 \check u_0^{N} + F_0 X_{-i}^{(N)}
+{\tfrac{1}{ N}}F_0 X_i^N\Big)dt+ D_0 d W_0(t) , \\
d  X_{-i}^{(N)}(t)& = \Big[\big(A+F\big)X_{-i}^{(N)} +B \check u_{-i}^{(N)} +\tfrac{1}{N}F  X_i^{N} +G  X_0^{N}\Big]dt\\
&\qquad \qquad + \tfrac{D}{N}\sum_{j\ne i} d W_{j}^N(t)
-\tfrac{1}{N}\Big(F X_{-i}^{(N)}+\tfrac{1}{N}F  X_i^{N} +G  X_0^{N}\Big)dt, \\
d  X_i^{N}(t)& = \Big( A  X_i^{N} +B  u_i^{N} +F  X_{-i}^{(N)}  +G  X_0^{N}\Big)dt  +D d W_i(t) +\tfrac{1}{N}FX_i^{N}dt,
\end{align*}
where  we use the notation $Y_{-i}^{(N)}= \frac{1}{N}\sum_{j\ne i}^N Y_j$.

Let $\tilde X_j^N$, $0\le j\le N$, denote the state variations caused by $\tilde u_i^N$.
The state variation of  $\mathcal{A}_j$, $1\le j\ne i\le N$, is
$$
d\tilde X_j^{N}(t) = \Big( A \tilde X_j^{N} +F \tilde X^{(N)}_{-i} +\tfrac{1}{N}F \tilde X_i^{N}
+G \tilde X_0^{N} \Big)dt,
$$
where $\tilde X_j^N(0)=0$.
 This implies $\tilde X_j^{N} =\tilde X_{j'}^{N}$ for all $1\le j, j'\ne i$.
Now,
\begin{align}
d\tilde X_0^{N}(t) &= \Big( A_0 \tilde X_0^{N} +F_0 \tilde X^{(N)}_{-i} +
 F_0 \tfrac{\tilde X_i^{N}}{N} \Big)dt,\nonumber \\
d\tilde X_{-i}^{(N)}(t) &= \Big[\big(A+F\big) \tilde X^{(N)}_{-i} +
F \tfrac{\tilde X_i^{N}}{ N}
+G \tilde X_0^{N} \Big]dt-\tfrac{1}{ N}\Big(F \tilde X^{(N)}_{-i}+
F \tfrac{\tilde X_i^{N}}{ N}
+G \tilde X_0^{N}\Big)dt,\nonumber  \\
d{\tilde X_i^{N}(t)} &=  ( A {\tilde X}_i^{N} + B{\tilde u_i^{N}})dt+\Big( F \tilde X^{(N)}_{-i} + F \tfrac{\tilde X_i^{N}}{ N}
+G \tilde X_0^{N} \Big) dt, \label{Xitut}
\end{align}
where $\tilde X_0^{N}(0)=\tilde X_{-i}^{(N)}(0) ={\tilde X_i^{N}(0)} =0$,
and $\tilde u_i^N  \in L^2_{\cal F}(0, T; \mathbb{R}^{n_1})$.

\begin{lem} \label{lemma:tilx}
There exists a constant $C$ independent of $N$ such that
$$
\sup_{0\le t\le T}E\Big(\big|\tilde X_0^{N}(t)\big|^2+\big|\tilde X_{-i}^{(N)}(t)\big|^2+
\tfrac{1}{ N^2}\big|\tilde X_i^{N}(t)\big|^2\Big)\le {C\over N^2}E\int_0^T|\tilde u_i^{N}(t)|^2dt.
$$
\end{lem}
{\it Proof.}
By solving the linear ODE of $(\tilde X_0^N,   \tilde X_{-i}^{(N)}, \tilde X_i^N/N)$, we first have a uniform bound estimate on the fundamental solution matrix on $[0,T]$ and next obtain the estimate
\begin{align}\label{xxxn}
\sup_{0\le t\le T}( |\tilde X_0^N(t)| +   |\tilde X_{-i}^{(N)}(t)| + |\tilde X_i^N(t)/N|)\le \frac{C}{N}\int_0^T |\tilde u_i^N(s)| ds.
\end{align}
The lemma follows by applying Schwarz inequality. \qed

\begin{rem}
It is  seen that  $(\tilde X_0^N, \tilde X_{-i}^{(N)})$ and  $\tilde X_i^N$ have two different scales.
\end{rem}

When the control changes from $(\check u_i^N, \check u_{-i}^N)$ to
$(\check u_i^N+\tilde u_i^N, \check u_{-i}^N)$,
the first variations of the costs have the following form
\begin{align*}
&\frac12 \delta J_0 =E\int_0^T \chi_0 (t)dt + E \chi_{0f},\quad
\frac{\lambda}{2N} \delta J_i= E\int_0^T \chi_i (t)dt+E \chi_{if},\quad\\
&\frac{\lambda}{2N} \sum_{1\le j\ne i}\delta J_j = E\int_0^T \chi_{-i}(t) dt+E \chi_{-if},
\end{align*}
where
\begin{align*}
 \chi_0=\ & \big[\check X_0^{N} -\big(H_0 \check X^{(N)} +\eta_0\big) \big]^T Q_0 \big(\tilde X_0^{N}-H_0
 \tilde X^{(N)}_{-i} -\mbox{$\frac{1}{N}$} H_0 \tilde X_i^{N}  \big ),\\
 \chi_i=\ &  \big[\check X_i^{N} -\big(H_1\check X_0^{N}+ H_2   \check X^{(N)} +\eta\big) \big]^T
\tfrac{1}{N} \lambda Q \Big(\tilde X_i^{N}- H_1 \tilde X_0^{N}-  H_2
 \tilde X^{(N)}_{-i} - \mbox{$\frac{1}{N}$}  H_2 \tilde X_i^{N}  \Big)\\
& + \big(\check u_i^{N}\big)^T \tfrac{1}{N} \lambda R
 \tilde u_i^{N},  \\
 \chi_{-i} = \ & \big[\big(I- H_2\big) \check X^{(N)} -H_1 \check X_0^{N} -\eta  \big]^T \lambda Q \Big[\big(I- H_2\big) \tilde X^{(N)}_{-i} -H_1 \tilde X_0^{N}  -\mbox{$\frac{1}{N}$}  H_2 \tilde X_i^{N}
  \Big]+ {\cal E}_1^{N}  ,
\end{align*}
and
\begin{align*}
 \chi_{0f}=\ & \big[\check X_0^{N} -\big(H_{0f} \check X^{(N)} +\eta_{0f}\big) \big]^T Q_{0f} \big(\tilde X_0^{N}-H_{0f}
 \tilde X^{(N)}_{-i} -\mbox{$\frac{1}{N}$} H_{0f} \tilde X_i^{N}  \big )(T),\\
 \chi_{if}=\ &  \big[\check X_i^{N} -\big(H_{1f}\check X_0^{N}+ H_{2f}   \check X^{(N)} +\eta_f\big) \big]^T
\tfrac{1}{N} \lambda Q_f \Big(\tilde X_i^{N}- H_{1f} \tilde X_0^{N}-  H_{2f}
 \tilde X^{(N)}_{-i} - \mbox{$\frac{1}{N}$}  H_{2f} \tilde X_i^{N}  \Big)(T),
  \\
 \chi_{-if} = \ & \big[\big(I- H_{2f}\big) \check X^{(N)} -H_{1f} \check X_0^{N} -\eta_f  \big]^T \lambda Q_f \Big[\big(I- H_{2f}\big) \tilde X^{(N)}_{-i} -H_{1f} \tilde X_0^{N}  -\mbox{$\frac{1}{N}$}  H_{2f} \tilde X_i^{N}
  \Big](T)+ {\cal E}^N_{f}.
\end{align*}
In the above,
\begin{align}
{\cal E}_1^N=& -\tfrac{\lambda}{N}(\check X_i^N )^TQ \Big[ (I-H_2) \tilde X_{-i}^{(N)} - H_1 \tilde X_0^{N} - \mbox{$\frac{1}{N}$}  H_2 \tilde X_i^{N}+
\tfrac{1}{N-1} \tilde X_{-i}^{(N)}  \Big] \nonumber \\
&+\tfrac{\lambda }{N-1}(\check X^{(N)})^T Q \tilde X_{-i}^{(N)}\nonumber  \\
&-\tfrac{\lambda }{N} ( H_1\check X_0^{N}+ H_2   \check X^{(N)} +\eta)^T
  Q \Big(  H_1 \tilde X_0^{N}+  H_2
 \tilde X^{(N)}_{-i} + \mbox{$\frac{1}{N}$}
 H_2 \tilde X_i^{N} \Big). \label{e1N}
\end{align}
See appendix B for the derivation of \eqref{e1N}.
The derivation of ${\cal E}^N_{f}$ is similar and omitted here.
We may regard ${\cal E}_1^N$ as  a higher order term relative to the first term in $\chi_{-i}$. Specifically, by Lemma \ref{lemma:tilx}, we have
$$
E|{\cal E}_1^N(t)|= O  \Big( \tfrac{1}{N^2} \Big (E [ |\check X_i^N(t)|^2
+ |\check X^{(N)}(t)|^2 +|\check X_0^N(t)|^2]\Big)^{1/2} \Big( E\int_0^T|\tilde u_i^{N}(t)|^2dt\Big)^{1/2} \Big).
$$
We may give a similar upper bound for $E|{\cal E}_f^N|$ by using $\big(E[|\check X_i^N(T)|^2+|\check X^{(N)}(T)|^2 +|\check X_0^N(T)|^2]\big)^{1/2}$ in place of the middle factor of $O(\cdot)$ above.

\begin{prop}\label{prop:pbp2} We have
\begin{align*}
E\int_0^T \big(\chi_0 +\chi_i +\chi_{-i}\big) dt + E( \chi_{0f} +\chi_{if} + \chi_{-if}) =0, \quad
\forall \tilde u_i^{N}  \in L^2_{\cal F}(0, T; \mathbb{R}^{n_1})  .
\end{align*}
\end{prop}

{\it Proof.} The proof is similar to that of Lemma \ref{lemma:pbp0} and we omit the detail. \qed

It can be shown that
\begin{align*}
 &\chi_0+\chi_i +\chi_{-i} \\
=\ & \Big[\check X_0^{N} -(H_0 \check X^{(N)} +\eta_0) \Big]^T Q_0 \Big(\tilde X_0^{N}-H_0
 \tilde X^{(N)}_{-i} -{\textstyle\frac{1}{N}} H_0 \tilde X_i^{N}  \Big) \\
 &  + \Big[\check X_i^{N} -(H_1\check X_0^{N}+ H_2   \check X^{(N)} +\eta) \Big]^T
\tfrac{1}{N} \lambda Q \tilde X_i^{N}
 + (\check u_i^{N})^T \tfrac{1}{N} \lambda R
 \tilde u_i^{N} \\
  &+ \Big[(I- H_2) \check X^{(N)} -H_1 \check X_0^{N} -\eta  \Big]^T \lambda Q
 \Big[(I- H_2) \tilde X^{(N)}_{-i} -H_1 \tilde X_0^{N}  -\tfrac{1}{N} H_2 \tilde X_i^{N}
  \Big]\\
  &+{\cal E}_2^{N},
\end{align*}
where ${\cal E}_2^N$ can again be treated as a higher order term,
and we may similarly rewrite $ \chi_{0f} +\chi_{if} + \chi_{-if}$.

\subsection{Limiting variational problem for the minor player}

Consider
\begin{align*}
&d\hat X_0^{\star}(t) = ( A_0 \hat X_0^{\star}+B_0 \hat u_0^{\star} +F_0 \hat m )dt
  + D_0 dW_0(t), \\
&d\hat m(t) = \big( (A+F) \hat m  +B \bar u
+G \hat X_0^{\star}\big)dt, \\
&d X_i^{\star}(t)= ( A X_i^{\star} + B \hat u_i^{\star} +F \hat m
+G\hat X_0^{\star} )dt +DdW_i(t),
\end{align*}
where $\hat u_0^\star$ has been determined from the solution of VP--(I), and $\hat u_i^\star\in L_{{\cal F}^{W_0,W_i}}^2(0,T; \mathbb{R}^{n_1})$.

Denote the state variational equations
\begin{align}\label{vastn}
\begin{cases}
d\tilde X_0^{\star}(t) = \Big( A_0 \tilde X_0^{\star} +F_0 \tilde m  +
\mbox{$\frac{1}{N}$}F_0 \tilde X_i^{\star} \Big)dt, \\
d\tilde m(t) = \Big[ (A+F) \tilde m  +\frac{1}{N}F \tilde X_i^{\star}
+G \tilde X_0^{\star} \Big]dt, \\
d\tilde X_i^{\star}(t) = ( A \tilde X_i^{\star} + B\tilde u_i^{\star}
)dt,
\end{cases}
\end{align}
where $\tilde X_0^{\star}(0)=\tilde m(0)= \tilde X_i^\star(0)=0$, and
$ \tilde u_i^\star\in L_{{\cal F}^{W_0,W_i}}^2(0,T; \mathbb{R}^{n_1}) $.

\begin{rem}
We see that  $(\tilde X_0^\star, \tilde m)$ and $ \tilde X_i^\star$ have different scales when $N$ increases,  which is similar to  the case of    $(\tilde X_0^N, \tilde X_{-i}^{(N)})$ and  $\tilde X_i^N$.
\end{rem}

We give some motivation for introducing the two variational equations in \eqref{vastn} containing the $1/N$ factor. For large $N$, if the perturbation $\tilde u_i^N$ is small, $J^{(N)}_{\rm soc}$ has a change by the order of $(1/N) (E\int_0^T |\tilde u_i^N|^2ds)^{\frac{1}{2}} $.
 Thus we need to look at the optimizing behavior at a finer scale. For this reason the two $1/N$ scaled terms in \eqref{vastn} are significant, and as it turns out below, they ensure that $(\tilde X_0^\star, \tilde m)$ provides a good approximation to $(\tilde X_0^N, \tilde X_{-i}^{(N)})$.

For approximating $\chi_0+\chi_i+\chi_{-i}$,  denote
\begin{align*}
\delta L_i^\star(t)=\ &\Big\{ [\hat X_0^{\star} -(H_0 \hat m +\eta_0) ]^T Q_0
\Big(\tilde X_0^{\star}-H_0
 \tilde m -\tfrac{1}{N} H_0 \tilde X_i^{\star}  \Big) \\
 &  + [\hat X_i^{\star} -(H_1\hat X_0^{\star}+ H_2   \hat m +\eta) ]^T
\tfrac{1}{N} \lambda Q \tilde X_i^{\star}
 + (\hat u_i^{\star})^T \tfrac{1}{N} \lambda R
 \tilde u_i^{\star} \\
  &+ [(I- H_2) \hat m -H_1 \hat X_0^{\star} -\eta  ]^T \lambda Q
   \Big[(I- H_2) \tilde m -H_1 \tilde X_0^{\star}  -\tfrac{1}{N}  H_2 \tilde X_i^{\star}
  \Big] \Big\} (t) \\
 =\ &  (\tilde X_0^\star)^ T (M_0\hat X_0^\star +K_0 \hat m+\nu_0) +\tilde m^T(K_0^T \hat X_0^\star+ M \hat m+\nu)\\
&+\tfrac{(\tilde X_i^\star)^T}{N}
(K_0^T \hat X_0^\star + (M-\lambda Q) \hat m+ \lambda Q \hat X_i^\star+ \nu)+ (\hat u_i^{\star})^T \tfrac{1}{N} \lambda R
 \tilde u_i^{\star} .
\end{align*}

In parallel to $\delta L_i^\star (t) $, we introduce a terminal variational term
\begin{align*}
\delta L_{if}^\star =  \ & \Big\{ (\tilde X_0^\star)^ T (M_{0f}\hat X_0^\star +K_{0f} \hat m+\nu_{0f}) +\tilde m^T(K_{0f}^T \hat X_0^\star+ M_f \hat m+\nu_f)\\
&+\tfrac{(\tilde X_i^\star)^T}{N}
(K_{0f}^T \hat X_0^\star + (M_f-\lambda Q_f) \hat m+ \lambda Q_f \hat
X_i^\star+\nu_f)\Big\} (T).
\end{align*}

{\bf  Variational Problem I\!I} VP--(I\!I):
Find $\hat u_i^{\star} \in L_{{\cal F}^{W_0,W_i}}^2(0,T; \mathbb{R}^{n_1})  $ such that
$$
E\int_0^T \delta L_i^\star dt +E \delta L_{if}^\star  =0, \qquad
 \forall\ \tilde u_i^{\star} \in L_{{\cal F}^{W_0,W_i}}^2(0,T; \mathbb{R}^{n_1})    .
 $$

The variational condition in VP--(I\!I) may be regarded as  an
approximation of the person-by-person optimality property as stated in Proposition \ref{prop:pbp2}.
The proposition below gives insights into the limiting variational problem VP--(I\!I) and provides a justification for the form of \eqref{vastn}.

\begin{prop}
 Let
$\tilde u_i^{N}= \tilde u_i^{\star} =v $ in \eqref{Xitut} and \eqref{vastn} for some fixed $v\in L^2_{\cal F}(0,T; \mathbb{R}^{n_1})$.
Then for some constant $C$ we have
\begin{align*}
\sup_{0\le t\le T}  E\Big[| \tilde X_0^{N}(t)-\tilde X_0^{\star}(t)|^2   + |\tilde X^{(N)}_{-i}(t)-\tilde m(t) |^2
 + |\tfrac{1}{N}\tilde X_i^{N}(t)-\tfrac{1}{N}\tilde X_i^{\star}(t)|^2\Big]\le {C\over N^4}.
\end{align*}
\end{prop}
{\it Proof.} Denote
$\delta_0(t)= \tilde X_0^N-\tilde X_0^\star $, $\delta_{-i}(t) = \tilde X_{-i}^{(N)}-\tilde m$, and $ \delta_i (t)= \tilde X_i^N-\tilde X_i^\star$.
Then we write
\begin{align*}
d\delta_0(t) &= \Big( A_0 \delta_0 +F_0  \delta_{-i}  +
 F_0 \tfrac{\delta_i}{ N} \Big)dt, \\
d\delta_{-i}(t)    &= \Big[\big(A+F\big) \delta_{-i}  + F
\tfrac{\delta_i }{ N}
+G\delta_0  \Big]dt-\tfrac{1}{N}\Big(F \tilde X^{(N)}_{-i}+ F
\tfrac{\tilde X_i^{N}}{ N}
+G \tilde X_0^{N}\Big)dt,\\
d \tfrac{\delta_i (t)}{N}  &=   A \tfrac{\delta_i}{N}  dt+ \tfrac{1}{N}\Big( F \tilde X^{(N)}_{-i} + F\tfrac{\tilde X_i^{N}}{ N}
+G \tilde X_0^{N} \Big) dt,
\end{align*}
where $ \delta_0(0)= \delta_{-i}(0) =\delta_i(0) =0$.
By assumption (A1),
\begin{align}
\sup_{0\le t\le T} |\delta_0 +\delta_{-i} +(\delta_i/N)|\le \frac{C}{N} \sup_{0\le t\le T}  ( |\tilde X^{(N)}_{-i}| + | \tilde X_i^{N}/ N|
+| \tilde X_0^{N}|).
\end{align}
Recalling \eqref{xxxn}, the proposition follows. \qed

For VP--(I\!I) and the associated variational equations in \eqref{vastn}, we try to identify  adjoint processes $(q_0, q, q_i)$ such that the equality in VP--(I\!I) is expressed only in terms of $\tilde u_i^\star$ and $(\hat X_0^\star, \hat m, \hat X_i^\star, \hat u_i^\star )$. Denote
\begin{align*}
dq_0(t)=\ & \psi_{11}dt + \psi _{12} dW_0(t)+\psi_{13} dW_i(t),  \\
dq (t)=\ & \psi_{21}dt + \psi _{22} dW_0(t)+\psi_{23} dW_i(t), \\
dq_i(t)=\ &\psi_{31}dt + \psi _{32} dW_0(t)+\psi_{33} dW_i(t),
\end{align*}
where $\psi_{jk}$ and $q_0(T)$, $q(T)$, $q_i(T)$ are to be determined.
After elementary although tedious computations, it turns out that $(q_0, q)$ and
$(p_0, p)$ in \eqref{clmaj} are determined by exactly  the same equations and terminal conditions. Thus, we may use the adjoint processes $(p_0, p, q_i)$ with the equation of $q_i$ appropriately determined.

Let $(\hat X_0^\star, \hat m, p_0, p)$ be solved first. After the above procedure of constructing the adjoint processes, we introduce the equation system
\begin{equation}\label{clMmin}
\begin{cases}
d{\hat X}_0^{\star}(t)= (A_0 {\hat X}_0^{\star}+B_0R_0^{-1}B_0^Tp_0 +F_0{\hat m})dt+D_0 dW_0(t), \\
d{\hat m}(t)= [(A+F) {\hat m} +B\bar u  +G {\hat X}_0^{\star}] dt,\\
d{\hat X}_i^{\star}(t)= ( A {\hat X}_i^{\star} + BR_\lambda^{-1}B^Tq_i +F{\hat m}
+G{\hat X}_0^{\star} )dt +DdW_i(t),\\
dp_0(t)  =(M_0\hat X_0^\star+K_0\hat m-A_0^Tp_0-G^Tp+\nu_0)dt+\xi_0(t)dW_0(t),\\
dp (t)  =[K_0^T\hat X_0^\star+M\hat m-F_0^Tp_0-(A+F)^Tp+ \nu] dt +\xi (t) dW_0(t),  \\
dq_i(t)  =\big[K_0^T{\hat X}_0^\star+(M-\lambda Q){\hat m}+\lambda Q{\hat X}_i^{\star} -F_0^Tp_0-F^Tp-A^Tq_i+\nu\big]dt\\
\qquad\qquad +\zeta_{ai}(t)dW_0(t)+\zeta_{bi}(t)dW_i(t),
\end{cases}
\end{equation}
where
\begin{align*}
&\hat X_0^\star(0)= X_0^N(0), \quad \hat m(0)= m_0,\quad \hat X_i^\star(0)= X_i^N(0) ,\\
&p_0(T)=- (M_{0f}\hat X_0^\star(T) +K_{0f} \hat m(T)+\nu_{0f}), \quad
 p(T)=-(K_{0f}^T \hat X_0^\star(T)+ M_f \hat m(T)+\nu_f),  \\
&q_i(T)=-(K_{0f}^T \hat X_0^\star(T) +
(M_f-\lambda Q_f)\hat m(T)+ \lambda Q_f \hat X_i^\star(T)+\nu_f) .
\end{align*}

\begin{thm} \label{theorem:vp2} Given  $\bar u\in L_{{\cal F}^{W_0}}^2(0,T; \mathbb{R}^{n_1})$,
 \eqref{clMmin} has a unique solution
$$(\hat X_0^\star, \hat m,  \hat X_i^\star, p_0, p, q_i, \xi_0, \xi, \zeta_{ai}, \zeta_{bi}  )$$ such that
\begin{align*}
&(\hat X_0^\star, \hat m,   p_0, p,  \xi_0, \xi) \in L_{{\cal F}^{W_0}}^2(0,T; \mathbb{R}^{4n})
 \times L_{{\cal F}^{W_0}}^2 (0,T;\mathbb{R}^{2n\times n_2}),\\
&(  \hat X_i^\star, q_i, \zeta_{ai}, \zeta_{bi}  ) \in L_{{\cal F}^{W_0,W_i}}^2(0,T; \mathbb{R}^{2n})
 \times L_{{\cal F}^{W_0,W_i}}^2 (0,T;\mathbb{R}^{2n\times n_2}),
\end{align*}
 and {\em VP-(I\!I)} has a unique solution given by
\begin{align}\label{vp2sol}
 \hat u_i^\star(t)=  R_\lambda^{-1}B^Tq_i(t).
 \end{align}
\end{thm}

{\it Proof.} Note that $(\hat X_0^\star, \hat m, p_0, p, \xi_0, \xi )$ is uniquely determined by Theorem \ref{theorem:4equ}.
To proceed, denote
\begin{align*}
\chi_1(t)&=G\hat X_0^\star+F\hat m,\\
\chi_2(t)&=K_0^T\hat X_0^\star+(M-\lambda Q)\hat m-F_0^Tp_0-F^Tp+\nu.
\end{align*}
We rewrite
\begin{align}
\label{eq2xiqi}
\begin{cases}
d{\hat X}_i^{\star}(t)= ( A {\hat X}_i^{\star} + BR_\lambda^{-1}B^Tq_i +\chi_1)dt +DdW_i(t),\\
dq_i(t)  =(\lambda Q{\hat X}_i^{\star}-A^Tq_i+\chi_2)dt +\zeta_{ai}(t)dW_0(t)+\zeta_{bi}(t)dW_i(t),
\end{cases}
\end{align}
for which we obtain a unique solution by using Lemma \ref{lemma:fbXYZ}
with the vector Brownian motion $(W_0, W_i)$.

We proceed to  show that \eqref{vp2sol} is a solution to VP--(I\!I), where the associated state processes are $\hat X_0^\star,\hat m,\hat X_i^\star$. Applying It\^o's formula to $ d[(\tilde X_0^\star)^T p_0 + \tilde m^T p +(\tilde X_i^\star/N)^T q_i] $ gives the relation
\begin{align}\label{dxmx}
E[(\tilde X_0^\star)^T p_0 + \tilde m^T p +(\tilde X_i^\star/N)^T q_i](T)= E\int_0^T \psi( \hat X_0^\star, \hat m, \hat X_i^\star, p_0, p, q_i )dt,
\end{align}
where the integrand may be easily determined.
Combining \eqref{dxmx} with \eqref{vp2sol}, we can show that $\hat u_i^\star$ satisfies the variational condition in VP--(I\!I).
The proof of uniqueness is similar to part ii) of Theorem \ref{theorem:4equ}. This is done by showing that a solution to VP--(I\!I) is necessarily represented as \eqref{vp2sol} via solving \eqref{clMmin}.  \qed

 To further analyze \eqref{eq2xiqi}, we introduce the backward stochastic Riccati differential equation (BSRDE)
\begin{align}
&-d{P}_\lambda(t)=(P_{\lambda}A+A^TP_{\lambda}
-P_{\lambda}BR_{\lambda}^{-1}B^TP_{\lambda}
+\lambda Q)dt -\sum_{k=1}^{n_2}
\Psi_k(t) dW_{0k}(t), \label{Plam} \\
&\qquad P_{\lambda}(T)=\lambda Q_f. \nonumber
\end{align}
By Lemma \ref{lemma:bsre},
we solve a unique $P_{\lambda}\ge0$  in
$L_{{\cal F}^{W_0}}^\infty(0,T; S^{n})$ with $\Psi_k \in L_{{\cal F}^{W_0}}^2(0,T; S^{n}) $.
Denote the BSDE
\begin{align*}
d\phi(t)=&\Big[\big(P_{\lambda}BR_{\lambda}^{-1}B^T-A^T\big)\phi
+P_{\lambda}\chi_1+\chi_2  \Big]dt+  \Lambda_0dW_0(t),
\end{align*}
where $\phi(T)=-(K_{0f}^T \hat X_0^\star(T) +
(M_f-\lambda Q_f)\hat m(T)+\nu_f)$.
We obtain a unique solution $(\phi, \Lambda_0)$ in $L_{{\cal F}^{W_0}}^2(0,T; \mathbb{R}^{n})
 \times L_{{\cal F}^{W_0}}^2 (0,T;\mathbb{R}^{n\times n_2})$.

\begin{lem}\label{lemma:zetaab}
 We have
$\zeta_{ai}=\Lambda_0- [\Psi_1\hat X_i^\star, \ldots,\Psi_{n_2}\hat X_i^\star]  $ and
$\zeta_{bi}=\zeta_b:=-P_\lambda D $.
\end{lem}

{\it Proof. }
 By the method in proving Lemma \ref{lemma:fbXYZ}, we can show
 $q_i$ is in fact given by $-P_\lambda \hat X_i^\star +\phi$. We further obtain    the relation
$$
P_\lambda D+\zeta_{bi}=0, \quad \zeta_{ai} =\Lambda_0- [\Psi_1\hat X_i^\star, \cdots,\Psi_{n_2}\hat X_i^\star]. 
$$
The lemma follows.
\qed

\section{Mean Field Social Optimum Solution}

\label{sec:mfsol}

\subsection{Consistency condition}
So far we have assumed that $\bar u(t)\in L_{{\cal F}^{W_0}}^2(0,T; \mathbb{R}^{n_1})  $, as the approximation of $\check u^{(N)}(t)$, is known for solving the variational problems VP--(I) and VP--(I\!I). Below we introduce a procedure to determine $\bar u$.
Let VP--(I\!I) be solved for $i=1, \ldots, N$, so that  \eqref{clMmin} determines
$$
\hat u_i^\star(t)=R_\lambda^{-1} B^T q_i (t),\qquad 1\le i\le N.
$$

Denote
$$q^{(N)}(t)={1\over N}\sum_{i=1}^Nq_i(t),\quad  \hat u^{\star (N)}(t) = \frac{1}{N} \sum_{i=1}^N \hat u_i^\star(t), \quad  \hat{ X}^{\star(N)}(t)={1\over N}\sum_{i=1}^N\hat{X}_i^\star(t).
$$
It is plausible to approximate $\bar u$ by
$\hat u^{\star (N)} (t)= R_\lambda^{-1} B^T q^{(N)} (t)$.
We obtain
\begin{align*}
dq^{(N)}(t)=&\big\{K_0^T\hat{ X}_0^\star+(M-\lambda Q)\hat{ m}+\lambda Q\hat{X}^{\star(N)}-F_0^Tp_0\\
&\qquad-F^Tp-A^Tq^{(N)}+\nu\big\}dt+\zeta_{a}^{(N)}dW_0(t)+\frac{1}{N}\sum_{i=1}^N
\zeta_bdW_i(t),
\end{align*}
where $\zeta_a^{(N)}=\frac{1}{N}\sum_{i=1}^N \zeta_{ai} $ and $\zeta_b$ is given in Lemma \ref{lemma:zetaab}.

Recall that $\hat m$ was introduced to  approximate $\hat X^{(N)}$.
Also, in view of Lemma \ref{lemma:zetaab}, let $\zeta_{a}^{(N)}(t)$ be approximated by $ \zeta_{a}(t)$.
When $N\to \infty$,  the above equation of $q^{(N)}$ is approximated by
\begin{align*}
d\bar q(t)=&\Big(K_0^T\hat{ X}_0^\star+M\hat{ m}-F_0^Tp_0-F^Tp-A^T\bar q+\nu\Big)dt+\zeta_{a}dW_0(t),
\end{align*}
where $\bar q(T)=-(K_{0f}^T \hat X_0^\star(T) +
M_f\hat m(T)+\nu_f)$. We solve a unique solution
$(\bar q, \zeta_a) \in L_{{\cal F}^{W_0}}^2(0,T; \mathbb{R}^{n})\times L_{{\cal F}^{W_0}}^2(0,T; \mathbb{R}^{n\times n_2})  $.
The proof of the next lemma is straightforward and omitted here.
\begin{lem}
We have $\bar q(t)=p(t)$ on $[0,T]$.
\end{lem}

Now we introduce  the following  consistency condition
\begin{align}
\bar u =  R_\lambda^{-1} B^T p. \label{consc}
\end{align}

We note that a fixed point property is embodied  in \eqref{consc}.
A similar situation also arises in mean field games \cite{HCM07,HMC06}.
Given a general $\bar u'\in  L_{{\cal F}^{W_0}}^2
(0,T;\mathbb{R}^{n_1}) $, we solve VP--(I)  to obtain a
well-defined adjoint process $p\in L_{{\cal F}^{W_0}}^2
(0,T;\mathbb{R}^{n_1})$, and  we  use an operator to denote
$\Gamma (\bar u') = R_\lambda^{-1} B^T p.$
So \eqref{consc} is equivalent to the fixed point relation
$\bar u= \Gamma (\bar u)$.

Typically in a mean field game with mixed players, one determines the consistency condition by combining the solutions of the two optimization problems of the major player and a representative minor player \cite{H10,NH12,NC13}. For the present problem, indeed we may determine $\bar u$ via $\bar q$ after solving VP--(I\!I). However, now $\bar p$ and $p$ coincide, and  for this reason \eqref{consc} is determined by the solution of VP--(I)  alone.

 \subsection{The system of FBSDEs}

Substituting $\bar u$ above into \eqref{clMmin}, we introduce the new system
\begin{equation}\label{eqcons}
\begin{cases}
d{\widehat X}_0^{\star}(t)= (A_0 {\widehat X}_0^{\star}+
F_0{\widehat m}+B_0R_0^{-1}B_0^Tp_0 )dt+D_0 dW_0(t), \\
d{\widehat m}(t)= [G {\widehat X}_0^{\star}+(A+F) {\widehat m}
+BR_\lambda^{-1}B^Tp] dt,\\
dp_0(t)  =(M_0{\widehat X}_0^\star+K_0{\widehat m}-A_0^Tp_0
-G^Tp+\nu_0)dt+\xi_0(t)dW_0(t),\\
dp (t)  =[K_0^T{\widehat X}_0^\star
+M{\widehat m}-F_0^Tp_0-(A+F)^Tp+\nu]dt +\xi (t) dW_0(t),
\end{cases}
\end{equation}
where $\widehat X_0^\star (0)= X_0^N(0)$, $\widehat m(0)=m_0$, $p_0(T)=-( M_{0f} \widehat X_0^\star(T) + K_{0f} \widehat m(T) + \nu_{0f})$,  $p(T)= -( K_{0f}^T \widehat X_0^\star(T) + M_f \widehat m(T) +\nu_f ) $; and its solution is used to define
\begin{align*}
\widehat \chi_1(t)&=G\widehat X_0^\star(t)+F\widehat m(t),\\
\widehat \chi_2(t)&=K_0^T\widehat X_0^\star(t)+(M-\lambda Q)\widehat m(t)-F_0^Tp_0(t)-F^Tp(t)+\nu.
\end{align*}

Note that \eqref{eqcons} differs from \eqref{clMmin} due to the elimination of $\bar u$ by the consistency condition. To distinguish the associated processes, we use the new notation $\widehat X_0^\star$ and $\widehat m$
in  \eqref{eqcons} in place of $\hat X_0^\star$ and $\hat m$.
However,
the variables $p_0,p$ are reused for the adjoint processes, and their identification should be clear from the context.

We further introduce
\begin{align} \label{eqcons2}
\begin{cases}
d{\widehat X}_i^{\star}(t)= \big[ A {\widehat X}_i^{\star}
+ BR_\lambda^{-1}B^Tq_i +\widehat \chi_1\big]dt +DdW_i(t),\\
dq_i(t)  =\big[\lambda Q{\widehat X}_i^{\star}-A^Tq_i+\widehat \chi_2\big]dt +\zeta_{ai}(t)dW_0(t)+\zeta_{b}(t)dW_i(t),
\end{cases}
\end{align}
where $\widehat X_i^\star(0)=X_i^N(0)$ and $q_i(T)=-(K_{0f}^T \widehat X_0^\star(T) +
(M_f-\lambda Q_f)\widehat m(T)+ \lambda Q_f \widehat X_i^\star(T)+\nu_f)$.

\begin{thm} \label{theorem:mainexist}
The  FBSDE \eqref{eqcons} has a unique solution
$
(\widehat X_0^\star, \widehat m, p_0, p, \xi_0, \xi)$
in
$$
 L_{{\cal F}^{W_0}}^2
(0,T;\mathbb{R}^{4n})   \times L_{{\cal F}^{W_0}}^2
(0,T;\mathbb{R}^{2n\times n_2}),
$$
and subsequently we uniquely solve \eqref{eqcons2} to obtain $({\widehat X}_i^{\star}, q_i, \zeta_{ai}, \zeta_b)$ in $$ L_{{\cal F}^{W_0,W_i}}^2 (0,T;\mathbb{R}^{2n}) \times L_{{\cal F}^{W_0,W_i}}^2 (0,T;\mathbb{R}^{2n\times n_2}).$$
\end{thm}

{\it Proof.}
We follow the notation in \eqref{not1}--\eqref{not3} and further
 denote
$$
{\bf X}_0=
\begin{bmatrix}
\widehat X_0^\star\\
\widehat m
\end{bmatrix},  \quad
{\mathbb B}_1=\left[
\begin{matrix}
B_0&0\\
0&B
\end{matrix}
\right],
\quad
{\mathbb R}_1
=
\left[
\begin{matrix}
R_0&0\\
0&R
\end{matrix}
\right].
$$
We  rewrite the system \eqref{eqcons} in the form
$$
\begin{cases}
d{\bf X}_0(t)=\big[{\mathbb A}_0{\bf X}_0(t)+{\mathbb B}_1{\mathbb R}_1^{-1}{\mathbb B}_1^T{\bf Y}_0(t)\big]dt+{\mathbb D}_0dW_0(t),\\
d{\bf Y}_0(t)=\big[{\mathbb Q}_0{\bf X}_0(t)-{\mathbb A}_0^T{\bf Y}_0(t)+{\bf v}_0\big]dt+{ Z}(t)dW_0(t),
\end{cases}
$$
where ${\bf Y}_0(T)=-{\mathbb Q}_{0f} {\bf X}_0(T) - {\bf v}_{0f}$.
By Lemma \ref{lemma:fbXYZ}, we uniquely solve $({\bf X}_0,{\bf Y}_0,Z)$, and subsequently \eqref{eqcons2}.
This completes the proof.~\qed

 Since ${\mathbb Q}_0\ge0$ and ${\mathbb Q}_{0f}\ge0$, let $({\bf P}\ge0, \Psi_1,\ldots, \Psi_{n_2})$ be the unique solution to the BSRDE
$$
-d{{\bf P}}(t)=[{\bf P}(t){\mathbb A}_0+{\mathbb A}_0^T{\bf P}(t)-{\bf P}(t){\mathbb B}_1{\mathbb R}_1^{-1}{\mathbb B}_1^T{\bf P}(t)+{\mathbb Q}_0]dt
-\sum_{k=1}^{n_2} \Psi_k(t) dW_{0k}(t) ,\quad {\bf P}_0(T)={\mathbb Q}_{0f}.
$$
 We further uniquely solve
\begin{align*}
d\varphi(t)=&\Big[-{\mathbb A}_0^T\varphi +{\bf P}(t){\mathbb B}_1{\mathbb R}_1^{-1}{\mathbb B}_1^T\varphi+ {\bf v}_0 +\sum_{k=1}^{n_2} \Phi_k {\mathbb D}_{0k}^{\rm col} \Big]dt+\Lambda dW_0(t)
\end{align*}
with the terminal condition $\varphi(T)=- {\bf v}_{0f}$.
Then we can write ${\bf Y}_0(t)=-{\bf P}(t){\bf X}_0(t)+\varphi(t)$.


For \eqref{eqcons2}, denote
$$
\widehat{ X}^{\star(N)}(t)={1\over N}\sum_{i=1}^N\widehat{ X}_i^{\star}(t),\quad q^{(N)}(t)={1\over N}\sum_{i=1}^Nq_i(t), \quad 0\le t\le T.
$$

\begin{lem}
\label{lemma:mfmq}
For \eqref{eqcons} and \eqref{eqcons2}, there exists a constant $C$ independent of $N$ such that
$$
\epsilon_{1,N}\mathrel{\mathop:}=
\sup_{0\le t\le T}E\big(\big|\widehat{X}^{\star(N)}(t)-\widehat{ m}(t)\big|^2
+\big|q^{(N)}(t)
-p(t)\big|^2\big)\le C (\tfrac{1}{ N}+|x_0^{(N)} -m_0|^2).
$$
\end{lem}
{\it Proof.}  Denote ${ y}_N(t)=\widehat{ X}^{\star(N)}(t)-\widehat{ m}(t)$ and ${ r}_N(t)=q^{(N)}(t)-p(t)$. We have
\begin{align*}
d{ y}_N(t)&=\Big[A{ y}_N(t)+BR_\lambda^{-1}B^T{ r}_N(t)\Big]dt+\frac{D}{N}
\sum_{i=1}^NdW_i(t),\\
d{ r}_N(t)&=\Big[\lambda Q { y}_N(t)-A^T{ r}_N(t)\Big]dt
+ (\zeta_a^{(N)}-\xi)dW_0(t)
+\frac{\zeta_b}{N}\sum_{i=1}^NdW_i(t),
\end{align*}
where $r_N(T)=q^{(N)}(T)-p(T)= \lambda Q_f (\widehat m(T)- \widehat X^{\star(N)}(T))=-\lambda Q_f y_N(T)$. Let $(P_\lambda, \Psi_1, \ldots, \Psi_{n_2})$ be solved from \eqref{Plam}.
Writing ${ r}_N(t)=-P_\lambda(t) { y}_N(t)+\psi_N(t)$, we obtain
\begin{align}
d\psi_N(t)&=\Big(P_\lambda(t)BR_\lambda^{-1}B^T-A^T\Big  )\psi_N(t)  dt
+ \sum_{k=1}^{n_2}[ \Psi_k y_N +(\zeta_a^{(N)}-\xi)_k^{\rm col}] dW_{0k}+
 \frac{P_\lambda D+\zeta_b}{N}
\sum_{i=1}^NdW_i(t)\nonumber\\
& = \Big(P_\lambda(t)BR_\lambda^{-1}B^T-A^T\Big  )\psi_N(t)  dt
+ \sum_{k=1}^{n_2}[ \Psi_k y_N +(\zeta_a^{(N)}-\xi)_k^{\rm col}] dW_{0k},\nonumber
\end{align}
where $ \psi_N(T)= q^{(N)}(T)-p(T) +\lambda Q y_N(T)=0$.
Note that $P_\lambda D+\zeta_b=0$ by Lemma \ref{lemma:zetaab}.  Take $Z$ with $Z_k^{\rm col} =\Psi_k y_N +(\zeta_a^{(N)}-\xi)_k^{\rm col}  $. Then $(\psi_N, Z)$ is a solution of the linear  BSDE.
This implies $\psi_N(t)=0$ and we further determine
$\Psi_k y_N +(\zeta_a^{(N)}-\xi)_k^{\rm col}=0$.
Next, by use of (A1) and ${ r}_N(t)=-P_\lambda(t) { y}_N(t)$,  we directly estimate $\sup_{0\le t\le T}E\big|{ y}_N(t)\big|^2$, which further gives a bound on $\sup_{0\le t\le T}E\big|{ r}_N(t)\big|^2$ since $P_\lambda$ is an  essentially bounded process. \qed

\section{Asymptotic Social Optimality}
\label{sec:asy}

Denote by ${\cal U}_{\rm centr}$  the set of centralized controls consisting of all
  $u=(u_0, u_1,\ldots, u_N)$, where each $u_j\in L^2_{\cal F}(0,T; \mathbb{R}^{n_1})$.  For a general $u\in {\cal U}_{\rm centr}$,  let the corresponding state processes be
  $(X_0^N, X_1^N,\ldots, X_N^N)$. We have the following equations
  \begin{align*}
dX_0^N(t) =& (A_0 X_0^N+B_0 u_0^N +F_0X^{(N)} )dt   +D_0 dW_0(t), \qquad \\
dX_i^N(t)=&  (A X_i^N +Bu_i^N+FX^{(N)}+GX_0^N ) dt
       +  DdW_i(t),
        \quad 1\leq i\leq N.
\end{align*}


We combine \eqref{eqcons} and \eqref{eqcons2} to write  the following FBSDE
\begin{equation}\label{eqn-consistency2}
\begin{cases}
d{\widehat X}_0^{\star}(t)= \big(A_0 {\widehat X}_0^{\star}+F_0{\widehat m}+B_0R_0^{-1}B_0^Tp_0 \big)dt+D_0 dW_0(t), \\
d{\widehat m}(t)= \big[G {\widehat X}_0^{\star}+(A+F) {\widehat m} +BR_\lambda^{-1}B^Tp\big] dt,\\
d{\widehat X}_i^{\star}(t)= \big( A {\widehat X}_i^{\star}+ BR_\lambda^{-1}B^Tq_i +F{\widehat m}
+G{\widehat X}_0^{\star} \big)dt +DdW_i(t),\\
dp_0(t)  =\big(M_0{\widehat X}_0^\star+K_0{\widehat m}-A_0^Tp_0-G^Tp
+\nu_0\big)dt
+\xi_0dW_0(t),\\
dp (t)  =\big[K_0^T{\widehat X}_0^\star+M{\widehat m}-F_0^Tp_0
-(A+F)^Tp+\nu\big]dt +\xi dW_0(t),\\
dq_i(t)  =\big[K_0^T{\widehat X}_0^\star+(M-\lambda Q){\widehat m}+\lambda Q{\widehat X}_i^{\star} -F_0^Tp_0-F^Tp-A^Tq_i+\nu\big]dt\\
\qquad\qquad +\zeta_{ai}(t)dW_0(t)+\zeta_{b}(t)dW_i(t),
\end{cases}
\end{equation}
where
 $\widehat X_0^\star (0)= X_0^N(0)$, $\widehat m(0)=m_0$,
 $\widehat X_i^\star (0)= X_i^N(0)$,
$p_0(T)=- (M_{0f}\widehat X_0^\star(T) +K_{0f} \widehat m(T)+\nu_{0f})$,
$ p(T)=-(K_{0f}^T \widehat X_0^\star(T)+ M_f \widehat m(T)+\nu_f)$,
$q_i(T)=-(K_{0f}^T \widehat X_0^\star(T) +
(M_f-\lambda Q_f)\widehat m(T)+ \lambda Q_f \widehat X_i^\star(T)+\nu_f) .$

We use Theorem \ref{theorem:mainexist} to determine the unique solution
 $({\widehat X}_0^{\star},\widehat{ m},\widehat{ X}_i^{\star},p_0,p ,q_i, \xi_0, \xi, \zeta_{ai}, \zeta_b\big)$ for \eqref{eqn-consistency2}.
Denote the set of individual controls
$$
\hat u_0^N=R_0^{-1}B_0^Tp_0, \quad
\hat u_i^N=R_\lambda^{-1}B^Tq_i, \quad 1\le i\le N.
$$
For $\hat u =\big(\hat u_0^N,\hat u_1^N,\ldots,\hat u_N^N\big)$,  let the corresponding
 state processes be $\big(\hat X_0^N, \hat X_1^N,\ldots,\hat X_N^N\big)$.
 \begin{align*}
d\hat X_0^N(t) =& \Big(A_0 \hat X_0^N+B_0 \hat u_0^N   +F_0\hat X^{(N)} \Big)dt   +D_0 dW_0(t), \qquad \\
d\hat X_i^N(t)=&  \Big(A \hat X_i^N +B\hat u_i^N+F\hat X^{(N)}+G\hat X_0^N \Big) dt
       +  DdW_i(t),
        \quad 1\leq i\leq N,
\end{align*}
where $\hat X_j^N(0)= X_j^N(0)$ for $0\le j\le N$.
It follows that
\begin{align}
d\hat X^{(N)}(t)=&  \Big[(A +F)\hat X^{(N)} +B\hat u^{(N)}+G\hat X_0^N \Big] dt
       + \frac{ D}{N}\sum_{i=1}^NdW_i(t).
\end{align}

 Denote $\tilde X^N_j(t)=X^N_j(t)-\hat X^N_j(t)$, $\tilde u^N_j(t)=u^N_j(t)-\hat u^N_j(t)$ for $0\le j\le N$, and
$$
 \tilde X^{(N)}=\frac{1}{N}\sum_{i=1}^N \tilde X_i^N,
\quad\tilde  u^{(N)}=\frac{1}{N}\sum_{i=1}^N \tilde u_i^N.
$$
We obtain
 \begin{align}
d\tilde X_0^N(t) =& \Big(A_0 \tilde X_0^N+B_0 \tilde u_0^N +F_0\tilde X^{(N)} \Big)dt, \label{erx0} \\
d\tilde X_i^N(t)=&  \Big(A \tilde X_i^N +B\tilde u_i^N+F\tilde X^{(N)}+G\tilde X_0^N \Big) dt, \quad 1\leq i\leq N,  \nonumber \\
d\tilde X^{(N)}(t)=&  \Big[\big(A+F\big) \tilde X^{(N)} +B\tilde u^{(N)}+G\tilde X_0^N \Big] dt,\label{erxN}
\end{align}
where $\tilde X_0^N(0)=  \tilde X_i^N(0) = \tilde X^{(N)}(0)=0$.

Denote
$$\epsilon_{2,N}= \sup_{0\le t\le T}E\Big(\big|\hat X^N_0(t)-\widehat X^\star_0(t)\big|^2+\big|\hat X^{(N)}(t)-\widehat m(t)\big|^2+\big|\hat u^{(N)}(t)-\bar u(t)\big|^2\Big),$$
where $\bar u=R_\lambda^{-1} B^T p$ and $p$ is given by \eqref{eqn-consistency2}.


\begin{lem} \label{lemma:en}
We have
$$
 \epsilon_{2,N}=O\Big(\tfrac{1}{ N}+| x_0^{(N)} -m_0|^2  \Big).
$$
\end{lem}

{\it Proof.} Note that $\hat u^{(N)}(t)-\bar u(t)=R_\lambda^{-1}B^T\big(q^{(N)}(t)-p(t)\big)$. Under (A1), Lemma \ref{lemma:mfmq} implies
$$
\sup_{0\le t\le T}E\big|\hat u^{(N)}(t)-\bar u(t)\big|^2=O\Big(\tfrac{1}{ N}+| x_0^{(N)} -m_0|^2   \Big).
$$
Denote
$$
y_0(t)=\hat X^N_0(t)-\widehat X^\star_0(t),\qquad y_N(t)=\hat X^{(N)}(t)-\widehat { X}^{\star(N)}(t).
$$
 Then $y_0(t)$ and  $y_N(t)$ satisfy the following linear ODE:
 $$
 {d\over dt}
 \left[
 \begin{matrix}
 y_0(t)\\
 y_N(t)
 \end{matrix}
 \right]
 =
  \left[
 \begin{matrix}
 A_0&F_0\\
G&A+F
 \end{matrix}
 \right]
  \left[
 \begin{matrix}
 y_0(t)\\
 y_N(t)
 \end{matrix}
 \right]
 +
  \left[
 \begin{matrix}
  F_0\big(\widehat{X}^{\star(N)}-\widehat m\big)\\
 F\big(\widehat{X}^{\star(N)}-\widehat m\big)
 \end{matrix}
 \right],\quad
 \left[
 \begin{matrix}
 y_0(0)\\
 y_N(0)
 \end{matrix}
 \right]
 =
 \left[
 \begin{matrix}
0\\
0
 \end{matrix}
 \right].
 $$
 Since all the parameter  processes are bounded and  $\sup_{0\le t\le T}E|\widehat{ X}^{\star(N)}(t)-\widehat m(t)|^2\le C(\tfrac{1}{ N}+| x_0^{(N)} -m_0|^2)$ by Lemma \ref{lemma:mfmq}, the lemma follows.
\qed

Now we are ready to state the asymptotic social optimality theorem.

\begin{thm} \label{theorem:tm}
We have
\begin{align}
\Big|J^{(N)}_{\rm soc}(\hat u)-\inf_{u\in {\cal U }_{\rm centr} }  J^{(N)}_{\rm soc}(u)\Big|=
O\Big(\tfrac{1}{\sqrt{N}}+ | x_0^{(N)} -m_0|  \Big). \nonumber
\end{align}
\end{thm}

The importance of the theorem comes from the fact that the set of decentralized individual controls $(\hat u_0, \hat u_1, \ldots, \hat u_N)$ can optimize $J_{\rm soc}^{(N)}(u)$ with little optimality loss in comparison with centralized controls.
The rest of this section is devoted to the proof of this theorem.

\subsection{Some Technical Lemmas}

Denote
\begin{align*}
&\Delta_0^N(t)=(\hat X_0^N-\Psi_0(\hat  X^{(N)}))^TQ_0
( \tilde X_0^N -H_0 \tilde X^{(N)}),\\
&\Delta_i^N(t)= ( \hat X_i^N - \Psi(\hat X_0^N, \hat X^{(N)}))^T Q (\tilde X_i^N- H_1\tilde X_{0}^N- H_2 \tilde X^{(N)}), \\
& \Delta_{0f}^N=\big\{(\hat X_0^N- H_{0f} \hat X^{(N)}-\eta_{0f})^T
Q_{0f} (\tilde X_0^N- H_{0f} \tilde X^{(N)})\big\}(T),  \\
&\Delta_{if}^N= \big\{(\hat X_i^N- H_{1f}\hat X_{0}^N- H_{2f} \hat X^{(N)} -\eta_f )^T Q_f (\tilde X_i^N- H_{1f}\tilde X_{0}^N-
 H_{2f} \tilde X^{(N)})\big\}(T).
\end{align*}

\begin{lem} \label{lemma:jhuju}
For any $u\in {\cal U }_{\rm centr}$, we have
\begin{align*}
J_{\rm soc}^{(N)}(u)\ge \ & J_{\rm soc}^{(N)} (\hat u)
+ 2E\int_0^T \Big[\Delta_0^N+\frac{\lambda}
{N}\sum_{i=1}^N\Delta_i^N
 +\big(\hat u_0^N\big)^T R_0 \tilde u_0^N
 +\frac{\lambda}{N}\sum_{i=1}^N
 \big(\hat u_i^N\big)^TR\tilde u_i^N  \Big] dt \\
 &+2 E\Big( \Delta_{0f}^N+ \frac{\lambda}{N}\sum_{i=1}^N\Delta_{if}^N  \Big).
\end{align*}
\end{lem}

{\it Proof.}  We check the integrands of $J_0$ and $J_i$ to obtain
\begin{align*}
&\big|X_0^N-\Psi_0(X^{(N)})\big|_{Q_0}^2+\big(u_0^N\big)^T R_0u_0^N \\
=\ &  \big|\hat X_0^N -\Psi_0(\hat X^{(N)}) +\tilde X_0^N -H_0 \tilde X^{(N)}\big|^2_{Q_0}
+\big|\hat u_0^N +\tilde u_0^N\big|_{R_0}^2\\
\ge \ &  \big|\hat X_0^N -\Psi_0(\hat X^{(N)})\big|_{Q_0}^2+\big(\hat u_0^N\big)^TR_0\hat u_0^N  + 2\Delta_0^N +2\big(\hat u_0^N\big)^TR_0 \tilde u_0^N
\end{align*}
and similarly,
\begin{align*}
 &\big|X_i^N-\Psi(X_0^N,X^{(N)})\big|_Q^2 + \big(u_i^N\big)^T R u_i^N\\
\ge \ &  \big|\hat X_i^N - \Psi(\hat X_0^N, \hat X^{(N)})\big|^2_Q+\big(\hat u_i^N\big)^T
R\hat u_i^N      +2\Delta_i^N +2\big(\hat u_i^N\big)^T R \tilde u_i^N.
\end{align*}
 We further check the terminal costs in $J_0$ and $J_i$ to obtain the estimate. \qed


We give a prior estimate on $\tilde X_0$ and $\tilde X^{(N)}$.
By elementary estimate we can show that there exists a constant $\hat C_0$
independent of $N$ such that
$$
J^{(N)}_{\rm soc}(\hat u_0^N, \hat u_1^N \ldots, \hat u_N^N)\le \hat C_0.
$$
For the estimate below it suffices to consider a set of individual controls $u=(u_0^N,\ldots,  u_N^N)\in {\cal U }_{\rm centr} $  such that
 \begin{align}
 J^{(N)}_{\rm soc}( u_0^N, u_1^N, \ldots,  u_N^N)\le \hat C_0. \label{juc0}
 \end{align}
 Denote all $u$ satisfying \eqref{juc0} by the set ${\cal U}_0$.

\begin{lem}\label{lemma:xtbound}
For all $u\in {\cal U}_0$, there exists $C_1$ such that
\begin{align}
\sup_{0\le t\le T}E\Big(\big|\tilde X_0^N(t)\big|^2+|\tilde X^{(N)}(t)\big|^2\Big) dt \le C_1.  \nonumber
\end{align}
\end{lem}

{\it Proof.}
By use of (A1)--(A3) and  direct SDE estimates for \eqref{eqn-consistency2} we can show that,
$$
\sup_{0\le j\le N} E\int_0^T \big|\hat u_j^N(t)\big|^2 dt \le C.
$$
Thus for $u\in {\cal U}_0$ , we have
\begin{align*}
E\int_0^T \big|\tilde u_j^N(t)\big|^2 dt &\le2 E\int_0^T \Big(\big|u_j^N(t)\big|^2 +\big|\hat u_j^N(t)\big|^2\Big)dt \le 2 E\int_0^T \big|u_j^N(t)\big|^2dt+C,\quad 0\le j\le N.
\end{align*}
Since $R_0(t),R(t)\ge c_1 I$ by (A1), \eqref{juc0} implies
$$
E\int_0^T\Big( \big|\tilde u_0^N(t)\big|^2+\big|\tilde u^{(N)}(t)\big|^2\Big)dt\le 2 E\int_0^T \Big(\big|u_0^N(t)\big|^2 +{1\over N}\sum_{i=1}^N\big| u_i^N(t)\big|^2\Big)dt
+C \le C_2,
$$
 where $C_2$ depends on $\hat C_0$.

By  \eqref{erx0} and \eqref{erxN}, we obtain for any $0\le t\le T$,
$$
|\tilde X_0^N(t)\big|+\big|\tilde X^{(N)}(t)\big|
  \le C \int_0^T\Big(\big|\tilde u_0^N(t)\big|+\big|\tilde u^{(N)}(t)\big|\Big) dt,
$$
and by applying Schwarz inequality,
\begin{align*}
E\Big(\big|\tilde X_0^N(t)\big|^2+\big|\tilde X^{(N)}(t)\big|^2\Big)
&  \le C E\int_0^T\Big(\big|\tilde u_0^N(t)\big|^2+\big|\tilde u^{(N)}(t)\big|^2\Big) dt
 \le C_1.
\end{align*}
This completes the proof.
\qed

Denote
\begin{align*}
\Theta(t)=\  &\Big\{ (\tilde X_0^N)^T (M_0\widehat X_0^\star +K_0 \widehat m +\nu_0) +\hat u_0^N R \tilde u_0^N \\
&+ (\tilde X^{(N)} )^T [K_0^T \widehat X_0^\star +(M-\lambda Q)\widehat m  +\nu ]
 + \frac{\lambda}{N}\sum_{i=1}^N [(\tilde X_i^N)^T Q \widehat X_i^\star + (\tilde u_i^N)^T R \hat u_i^N]\Big\}(t),
\end{align*}
and
\begin{align*}
\Theta_f= \ &\Big\{ (\tilde X_0^N)^T (M_{0f}\widehat X_0^\star +K_{0f}
\widehat m +\nu_{0f}) \\
&+ (\tilde X^{(N)} )^T [K_{0f}^T \widehat X_0^\star +(M_f-\lambda Q_f)\widehat m  +\nu_f  ]
 + \frac{\lambda}{N}\sum_{i=1}^N (\tilde X_i^N)^T Q_f
 \widehat X_i^\star\Big\}(T) .
\end{align*}
\begin{lem} \label{lemma:Theta}
Suppose $u\in {\cal U }_{\rm centr}$.  Then
\begin{align} \label{ththf}
E\int_0^T \Theta (t) dt+E \Theta_f +E\int_0^T [ (\tilde X^{(N)} )^T F^T+(\tilde X_0^N)^TG^T ](q^{(N)}-p) dt  =0.
\end{align}
If, in addition, $u\in {\cal U}_0$, then
$$
\Big|E\int_0^T \Theta (t) dt +E\Theta_f \Big| \le C \left(E \int_0^T
|q^{(N)}-p|^2dt\right)^\frac12.
$$
\end{lem}

{\it Proof.}
We have
\begin{align*}
d(p_0^T \tilde X_0^N)=&\Big[(\tilde X_0^N)^T (M_0\widehat X_0^\star +K_0
\widehat m -G^Tp+\nu_0) +p_0^T(B_0 \tilde u_0^N+F_0\tilde X^{(N)}) \Big]dt \\
&+ (\tilde X_0^N)^T \xi_0dW_0.  \end{align*}
Then
\begin{align} \label{eqg1}
E\int_0^T g_{1}(t)dt+ E\{(\tilde X_0^N(T))^T(M_{0f}\widehat X_0^\star(T) +K_{0f}
\widehat m(T) +\nu_{0f})\} =0,
\end{align}
where
$$
g_1(t)=(\tilde X_0^N)^T (M_0\widehat X_0^\star +K_0 \widehat m -G^Tp+\nu_0) +p_0^T(B_0 \tilde u_0^N+F_0\tilde X^{(N)}) .
$$
By checking $d  (\frac{1}{N}\sum_{i=1}^N q_i^T \tilde X_i^N)$, we obtain
  \begin{align}\label{eqg2}
E\int_0^T g_2 dt +\frac{1}{N}\sum_{i=1}^N E\Big\{(\tilde X_i^N(T))^T \big[ K_{0f}^T \widehat X_0^\star(T) +(M_f-\lambda Q_f)\widehat m (T)+\lambda Q_f \widehat X_i^\star(T)+\nu_f \big]\Big\}=0,
  \end{align}
where
\begin{align*}
g_2 (t)= & (\tilde X^{(N)} )^T [K_0^T \widehat X_0^\star +(M-\lambda Q)\widehat m-F_0^Tp_0 -F^T p +\nu +F^T q^{(N)} ] \\
& +  (\tilde X_0^N)^T G^T q^{(N)} + \frac{\lambda}{N}\sum_{i=1}^N (\tilde X_i^N)^T Q \widehat X_i^\star +\frac{1}{N} \sum_{i=1}^N (\tilde u_i^N)^T B^T q_i.
\end{align*}

Then
\begin{align*}
g_1(t)+g_2(t)=\ & [(\tilde X_0^N)^T (M_0\widehat X_0^\star +K_0 \widehat m +\nu_0) +\hat u_0^N R \tilde u_0^N ]\\
&+ (\tilde X^{(N)} )^T [K_0^T \widehat X_0^\star +(M-\lambda Q)\widehat m  +\nu  ]   \\
& + \frac{\lambda}{N}\sum_{i=1}^N (\tilde X_i^N)^T Q \widehat X_i^\star +\frac{\lambda}{N} \sum_{i=1}^N (\tilde u_i^N)^T R \hat u_i^N\\
& + (\tilde X^{(N)} )^T F^T(q^{(N)}-p)+(\tilde X_0^N)^T
G^T (q^{(N)}-p)\\
=\ &\Theta  + (\tilde X^{(N)} )^T F^T(q^{(N)}-p)+(\tilde X_0^N)^T
G^T (q^{(N)}-p).
\end{align*}
By
\eqref{eqg1}--\eqref{eqg2}, we derive \eqref{ththf}.
The remaining part follows by applying Schwarz theorem and Lemma
\ref{lemma:xtbound}. \qed

\subsection{Proof of Theorem \ref{theorem:tm}}

We have
\begin{align*}
&\Delta_0^N+ \frac{\lambda}{N}\sum_{i=1}^N\Delta_i^N\\
=\ & (\tilde X_0^N)^T \Big[( Q_0+\lambda H_1^T QH_1) \hat X_0^N  -(Q_0H_0 +\lambda H_1^TQ (I-H_2)) \hat X^{(N)}   +\lambda H_1^TQ\eta -Q_0\eta_0 \Big]\\
&+ (\tilde X^{(N)})^T [  H_0^T Q_0 H_0-\lambda Q H_2-\lambda H_2^T Q (I-H_2)] \hat X^{(N)}\\
&+ (\tilde X^{(N)})^T (\lambda H_2^T QH_1 -H_0^T Q_0 -\lambda QH_1) \hat X_0^N \\
&+(\tilde X^{(N)})^T (\lambda H_2^T Q\eta +H_0^T Q_0\eta_0 -\lambda Q\eta)\\
&+ \frac{\lambda}{N}\sum_{i=1}^N (\hat X_i^N)^T Q \tilde X_i^N\\
=\ &(\tilde X_0^N)^T (M_0 \hat X_0^N  +K_0 \hat X^{(N)}   + \nu_0)
+ (\tilde X^{(N)})^T [ (M-\lambda Q)\hat X^{(N)} +K_0^T \hat X_0^N + \nu   ]\\
  &+ \frac{\lambda}{N}\sum_{i=1}^N (\hat X_i^N)^T Q \tilde X_i^N.
 \end{align*}
Since
$$
d\big(\hat X^N_i(t)-\widehat X_i^\star(t)\big)=\Big[A(\hat X^N_i-\widehat X_i^\star)+F(\hat X^{(N)}-\widehat m)+G(\hat X^N_0-\widehat X_0^\star)\Big]dt,
$$
and $\hat X^N_i(0)-\widehat X_i^\star(0)=0$ for $1\le i\le N$,  we obtain
$$
\hat X^N_i(t)-\widehat X_i^\star(t)=\hat X^{(N)}(t)-\widehat{ X}^{\star(N)}(t).
$$

So
\begin{align}
&\Delta_0^N+\frac{\lambda}
{N}\sum_{i=1}^N\Delta_i^N
 +\big(\hat u_0^N\big)^T R_0 \tilde u_0^N
 +\frac{\lambda}{N}\sum_{i=1}^N
 \big(\hat u_i^N\big)^TR\tilde u_i^N -\Theta \nonumber \\
 =\ &(\tilde X_0^N)^T [M_0 (\hat X_0^N- \widehat X_0^\star  )
 +K_0 (\hat X^{(N)}-\widehat m) ] \nonumber \\
 &+(\tilde X^{(N)} )^T [K_0^T(\hat X_0^N- \widehat X_0^\star  ) +(M-\lambda Q)(\hat X^{(N)}-\widehat m)  ] \nonumber \\
&+\frac{\lambda}{N}\sum_{i=1}^N (\tilde X_i^N)^T Q  (\hat X_i^N -\widehat X_i^\star) \nonumber  \\
=\ & (\tilde X_0^N)^T [M_0 (\hat X_0^N- \widehat X_0^\star  )
 +K_0 (\hat X^{(N)}-\widehat m) ] \nonumber \\
 &+(\tilde X^{(N)} )^T [K_0^T(\hat X_0^N- \widehat X_0^\star  ) +M(\hat X^{(N)}-\widehat m) -\lambda Q (\widehat X^{\star(N)}-\widehat m)]. \label{dth}
\end{align}
In a similar manner, we can show
\begin{align}
\Delta_{0f}^N+ \frac{\lambda}{N}\sum_{i=1}^N\Delta_{if}^N -\Theta_f=\ &\Big\{ (\tilde X_0^N)^T [M_{0f}(\hat X_0^N- \widehat X_0^\star  )
 +K_{0f} (\hat X^{(N)}-\widehat m) ] \label{dthf} \\
 &+(\tilde X^{(N)} )^T [K_{0f}^T(\hat X_0^N- \widehat X_0^\star  ) +M_f(\hat X^{(N)}-\widehat m) -\lambda Q_f (\widehat X^{\star(N)}-\widehat m)]\Big\}(T). \nonumber
\end{align}

 It follows from   \eqref{dth}--\eqref{dthf} that
\begin{align*}
K_1\mathrel{\mathop:}=&\Big|E\int_0^T\Big[\Delta_0^N+\frac{\lambda}
{N}\sum_{i=1}^N\Delta_i^N
 +\big(\hat u_0^N\big)^T R_0 \tilde u_0^N
 +\frac{\lambda}{N}\sum_{i=1}^N
 \big(\hat u_i^N\big)^TR\tilde u_i^N\Big] dt+E(\Delta_{0f}^N+ \frac{\lambda}{N}
 \sum_{i=1}^N\Delta_{if}^N)\Big|\\
 \le\ & \Big|E\int_0^T \Theta(t) dt +E\Theta_f\Big|
  + CE\int_0^T \phi_N(t) dt  +C E\phi_N(T),
 \end{align*}
where $\phi_N(t)=\{\big(|\tilde X_0^N|+ | \tilde X^{(N)}|\big) \big(|\hat X_0^N -\widehat X_0^\star|+|\hat X^{(N)}-\widehat m|
+|\widehat X^{\star(N)}-\widehat m| \big)\}(t)$.
Lemmas \ref{lemma:xtbound} and \ref{lemma:Theta} imply that
\begin{align}
K_1\le C (\epsilon_{1,N}+\epsilon_{2,N})^{1/2}. \nonumber
\end{align}

By Lemma \ref{lemma:jhuju} and the above upper bound for $K_1$, for all $u\in {\cal U}_0$, we have
\begin{align}
J^{(N)}_{\rm soc}(\hat u) \le J^{(N)}_{\rm soc}(u)   +O((\epsilon_{1,N}+\epsilon_{2,N})^{1/2}),
\label{jjep}
\end{align}
which is automatically true when $u$ is not in ${\cal U}_0$.
Recalling Lemmas \ref{lemma:mfmq} and  \ref{lemma:en}, we complete the proof.
 \qed

\section{Conclusion}
\label{sec:con}

This paper studies an LQ mean field social optimization problem with mixed players. The solution is obtained by exploiting a person-by-person optimality principle and constructing two low dimensional limiting variational problems. This method derives an FBSDE system for the major player and a representative minor player. We prove the existence and uniqueness of the solution to the FBSDE and establish asymptotic social optimality for the resulting decentralized controls of the $N+1$ players.

\section*{Appendix A}

\renewcommand{\theequation}{A.\arabic{equation}}
\setcounter{equation}{0}
\renewcommand{\thethm}{A.\arabic{thm}}
\setcounter{thm}{0}

\begin{lem} \label{lemma:bsre}  {\em \cite{P92,KT03}}
Assume

 i) $\{\hat W(t)=[\hat W_1(t),\ldots, \hat W_l(t)]^T, t\ge 0\}$ is an $\mathbb{R}^l$-valued standard Brownian motion;

ii)  $\{\hat A(t), \hat B(t), \hat Q(t), \hat R(t), 0\le t \le T\}$ are ${\cal F}_t^{\hat W}$-adapted essentially bounded processes and are $\mathbb{R}^{k\times k}$, $\mathbb{R}^{k\times k_1}$, $S_+^{k}$, ${S}_+^{ k_1}$-valued, respectively;
 $\hat  R(t)\ge \alpha I$ for a deterministic constant $\alpha>0$; and $\hat Q_f$  is $S_+^k$-valued, ${\cal  F}_T^{\hat W}$-measurable, and  essentially bounded.


Then the backward stochastic Riccati differential  equation (BSRDE)
$$
\begin{cases}
-dP(t)=\big( \hat A^T P +P \hat A -P \hat B\hat R^{-1}\hat B^TP+\hat Q\big)(t)dt- \sum_{i=1}^l \Psi_i(t) d\hat W_i(t),\\
P(T)= \hat Q_f
\end{cases}
$$
has a unique ${\cal F}^{\hat W}_t$-adapted solution
$(P, \Psi_1, \ldots, \Psi_l)$  satisfying that  $P$ is $S^k_+$-valued and essentially bounded, and that each $\Psi_i\in  L^2_{\cal F^{\hat W}}(0,T; S^k)$.
\end{lem}

More general forms of this Riccati equation were studied in \cite[sec. 5]{P92}, \cite[sec. 2]{KT03},
where $\Psi_i$ also appears linearly in the drift term. The proof method was presented in \cite[sec. 5]{P92} by applying quasi-linearization of the Riccati equation.

We further introduce the assumption
\begin{align}\label{gvdv}
g,v\in L^2_{\cal F^{\hat W}}(0,T;
 {\mathbb R}^k) , \quad \hat D \in L^2_{\cal F^{\hat W}}(0,T;
 {\mathbb R}^{k\times l}), \quad v_f\ \mbox{is}\ {\cal F}_T^{\hat W}\mbox{-measurable},\ E|v_f|^2<\infty.
\end{align}

Consider the FBSDE
\begin{align}\label{fbxyh}
\begin{cases}
dX(t)=(\hat AX+\hat B\hat R^{-1}\hat B Y +g) dt +\hat D d\hat W(t),\\
dY(t)=( \hat QX -\hat A^T Y +v ) dt +Z d\hat W(t),
\end{cases}
\end{align}
where $Y(T)=-\hat Q_f X(T) -v_f$ and $X(0)=x_0\in \mathbb{R}^k$.

Denote the linear BSDE
\begin{align*}
d\psi(t)=&\big(-{\hat A}^T\psi +{ P}\hat{ B}{\hat R}^{-1}{\hat B}^T\psi+ { P}g  +{ v} +\sum_{i=1}^l
\Psi_i \hat D_{i}^{\rm col}\big)dt+ \Lambda d\hat W(t),
\end{align*}
where $\psi(T)= - { v}_{f}$.
There exists a unique solution $(\psi, \Lambda )\in L^2_{{\cal F}^{\hat W}}(0,T;
 {\mathbb R}^k)\times
 L^2_{{\cal F}^{\hat W}}(0,T; {\mathbb R}^{k\times l})$.

\begin{lem}\label{lemma:fbXYZ}
Suppose the assumptions in Lemma \ref{lemma:bsre}
and \eqref{gvdv} hold,
then \eqref{fbxyh} has a unique solution $(X,Y,Z)$ in $ L^2_{{\cal F}^{\hat W}}(0,T;
 {\mathbb R}^{2k})\times L^2_{{\cal F}^{\hat W}}(0,T;
 {\mathbb R}^{k\times l})$,
and
$$
Y= -PX +\psi, \qquad Z_i^{\rm col}=\Lambda_i^{\rm col} -P\hat D_i^{\rm col} -\Psi_iX.
$$
\end{lem}

{\it Proof.}
To show existence, consider the SDE
$$
dX(t)=(\hat AX+\hat B\hat R^{-1}\hat B (-PX+\psi) +g) dt +\hat D d\hat W(t),
\quad X(0)=x_0,
$$
which has  a unique solution; we  choose $Y=-PX+\psi$.
By It\^o's formula, we derive
\begin{align}
dY=(\hat QX-\hat A^T Y +v) dt +\sum_{i=1}^l [\Lambda_i^{\rm col} -P\hat D_i^{\rm col} -\Psi_iX] d\hat W_i.
\end{align}
We choose $Z_i^{\rm col}=[\Lambda_i^{\rm col} -P\hat D_i^{\rm col} -\Psi_iX]  $ for all $i\le l$. Then $(X,Y, Z)$ constructed above is a
  solution to \eqref{fbxyh}.

To show uniqueness,
suppose there is another solution $(X',Y',Z')$.
Denote
$\tilde X= X-X'$, $\tilde Y= Y-Y'$  and $\tilde Z= Z-Z'$. So $\tilde Y(T)= -\hat Q_f\tilde  X(T)$.
Denote $\tilde Y= -P \tilde X +\tilde \varphi$, where $\tilde \varphi$ is to be determined.
By It\^o's formula,
$$
d\tilde \varphi = (-\hat A^T \tilde \varphi +P\hat B \hat R^{-1} \hat B^T \tilde\varphi) dt
+ \tilde Z d\hat W +\sum_{i=1}^l \Psi_i \tilde X d\hat W_i,
$$
where $\tilde \varphi (T)=0$. Note that $\tilde X$ has been given and $(\tilde \varphi, \tilde Z)$ is a solution to the above linear BSDE. We necessarily have $\tilde \varphi =0$ and
$\tilde Z_i^{\rm col}=-\Psi_i \tilde X$.
We can further show $\tilde X=\tilde Y=0$ and $\tilde Z=0$. This proves uniqueness.
\qed

\section*{Appendix B}
{\it Derivation of \eqref{e1N}}:
We have the first order cost variation: For $1\le j\ne i$,
$$
\tfrac{1}{2}\delta J_j = \big[\check X_j^{N} -\big(H_1\check X_0^{N}+ H_2   \check X^{(N)} +\eta\big) \big]^T
\tfrac{1}{N}  Q \Big(\tilde X_j^{N}- H_1 \tilde X_0^{N}-  H_2
 \tilde X^{(N)}_{-i} - \mbox{$\frac{1}{N}$}  H_2 \tilde X_i^{N}  \Big).
$$
 By the fact that all $\tilde X_j$, $1\le j\ne i$, are equal,
we calculate
\begin{align*}
\Delta_1&\mathrel{\mathop:} =\sum_{1\le j\ne i} \big(\check X_j^{N}  \big)^T
\tfrac{1}{N}\lambda  Q \Big(\tilde X_j^{N}- H_1 \tilde X_0^{N}-  H_2
 \tilde X^{(N)}_{-i} - \mbox{$\frac{1}{N}$}  H_2 \tilde X_i^{N}  \Big)\\
 &=\big(\check X_{-i}^{(N)}  \big)^T  \lambda Q \Big(\tilde X_j^{N}- H_1 \tilde X_0^{N}-  H_2
 \tilde X^{(N)}_{-i} - \mbox{$\frac{1}{N}$}  H_2 \tilde X_i^{N}  \Big)\\
 &= \big(\check X^{(N)} -\tfrac{1}{N} \check X_i^N \big)^T  \lambda Q
 \Big(\tfrac{N}{N-1} \tilde X_{-i}^{(N)}- H_1 \tilde X_0^{N}-  H_2
 \tilde X^{(N)}_{-i} - \mbox{$\frac{1}{N}$}  H_2 \tilde X_i^{N}  \Big)\\
 &= \big(\check X^{(N)} -\tfrac{1}{N} \check X_i^N \big)^T  \lambda Q
 \Big( (I-H_2) \tilde X_{-i}^{(N)} - H_1 \tilde X_0^{N} - \mbox{$\frac{1}{N}$}  H_2 \tilde X_i^{N}+ \tfrac{1}{N-1} \tilde X_{-i}^{(N)} \Big) ,
\end{align*}
and
\begin{align*}
\Delta_2&\mathrel{\mathop:} =\sum_{1\le j\ne i} \big( H_1\check X_0^{N}+ H_2   \check X^{(N)} +\eta\big)^T
\tfrac{\lambda}{N}  Q \Big(\tilde X_j^{N}- H_1 \tilde X_0^{N}-  H_2
 \tilde X^{(N)}_{-i} - \mbox{$\frac{1}{N}$}  H_2 \tilde X_i^{N}  \Big)\\
 &=\big( H_1\check X_0^{N}+ H_2   \check X^{(N)} +\eta\big)^T
 \\
 &\quad \lambda Q \Big(\tilde X_{-i}^{(N)} - H_1 \tilde X_0^{N}-  H_2
 \tilde X^{(N)}_{-i} - \mbox{$\frac{1}{N}$}  H_2 \tilde X_i^{N} +\tfrac{1}{N}(H_1 \tilde X_0^{N}+  H_2
 \tilde X^{(N)}_{-i} + \mbox{$\frac{1}{N}$}  H_2 \tilde X_i^{N}) \Big)\\
 &=\big( H_1\check X_0^{N}+ H_2   \check X^{(N)} +\eta\big)^T
 \\
 &\quad \lambda Q \Big( (I-H_2) \tilde X_{-i}^{(N)} - H_1 \tilde X_0^{N}
- \mbox{$\frac{1}{N}$}  H_2 \tilde X_i^{N} +\tfrac{1}{N}(H_1 \tilde X_0^{N}+  H_2
 \tilde X^{(N)}_{-i} + \mbox{$\frac{1}{N}$}  H_2 \tilde X_i^{N}) \Big).
\end{align*}

We may write
\begin{align*}
\frac{\lambda}{2N}\sum_{j\ne i}^N\delta J_j =\Delta_{1} -\Delta_2.
\end{align*}
Subsequently, we  determine the form of ${\cal E}_1^N$ as in \eqref{e1N}.

\end{document}